\documentclass[10pt,twoside]{amsart}
\usepackage{amsmath,amssymb,amsthm,amscd,latexsym,graphicx}
\usepackage{epigraph}
\usepackage[OT2,OT1]{fontenc}
\usepackage{graphicx}
\usepackage[all]{xy}
\usepackage{hyperref}
\usepackage{enumerate}

\usepackage{xcolor}

 \font\caps=cmcsc10                    
 \font\Caps=cmcsc10 scaled \magstep1   

 \makeatletter
 \@specialpagefalse\headheight=8.5pt
 \makeatother

\def\TSkip{\medskip}
 \newbox\TheTitle{\obeylines\gdef\GetTitle #1
 \ShortTitle  #2
 \SubTitle    #3
 \Author      #4
 \ShortAuthor #5
 \EndTitle
 {\setbox\TheTitle=\vbox{\baselineskip=20pt\let\par=\cr\obeylines%
 \halign{\centerline{\Caps##}\cr\noalign{\medskip}\cr#1\cr}}%
         \copy\TheTitle\TSkip\TSkip%
 \def\next{#2}\ifx\next\empty\gdef\STitle{#1}\else\gdef\STitle{#2}\fi%
 \def\next{#3}\ifx\next\empty%

 \else\setbox\TheTitle=\vbox{\baselineskip=20pt\let\par=\cr\obeylines%
     \halign{\centerline{\caps##} #3\cr}}\copy\TheTitle\TSkip\TSkip\fi%
 \centerline{\caps #4}\TSkip\TSkip%
 \def\next{#5}\ifx\next\empty\gdef\SAuthor{#4}\else\gdef\SAuthor{#5}\fi%
 \catcode'015=5}}

 \def\Abstract{\begingroup\narrower
     \parskip=\medskipamount\parindent=0pt{\caps Abstract. }}

 \long\def\MSC#1\EndMSC{\def\arg{#1}\ifx\arg\empty\relax\else
      {\par\narrower\noindent%
      2000 Mathematics Subject Classification: #1\par}\fi}

 \long\def\KEY#1\EndKEY{\def\arg{#1}\ifx\arg\empty\relax\else
         {\par\narrower\noindent Keywords and Phrases: #1\par}\fi\TSkip}

 \long\def\DATE#1\EndDATE{\def\arg{#1}\ifx\arg\empty\relax\else
         {\par\narrower\noindent \center{\textit{#1}}\par}\fi\TSkip\TSkip\TSkip}

 \font\bf= cmbx10 at 10pt

\newcommand{\Sym}{\hbox{\ptitgot S}}

\renewcommand{\to}{\longrightarrow}

\newcommand{\mylabel}[2]{#2\def\@currentlabel{#2}\label{#1}}

 \newcommand{\F}{\mathbb{F}}
 \newcommand{\G}{\mathbb{G}}
 \newcommand{\Gm}{\mathbb{G}_{\mathrm {m}}}
  \newcommand{\Ga}{\mathbb{G}_{\mathrm {a}}}

 \renewcommand{\P}{\mathbb P}
 \newcommand{\Q}{\mathbb{Q}}
 
 \newcommand{\Z}{\mathbb{Z}}

 \newcommand{\GL}{\mathrm{GL}}
 
 \renewcommand{\O}{\mathcal O}

  \newcommand{\Det}{\mathrm{Det}}

 \newcommand{\Aut}{\mathrm{Aut}}

 \newcommand{\Fl}{\mathrm{Fl}}
 \newcommand{\carac}{\mathrm{char}}

  \newcommand{\End}{\mathrm{End}}
 \newcommand{\Ext}{\mathrm{Ext}}

 \newcommand{\Hom}{\mathrm{Hom}}
 
  \newcommand{\Ver}{\mathrm{Ver}}
   \newcommand{\Frob}{\mathrm{Frob}}

 \newcommand{\M}{\mathcal{M}}

\newcommand{\Id}{\mathrm{Id}}

\newcommand{\EExal}{\mathbf{Exal}}

\renewcommand{\Im}{\mathrm{Im}}

\newcommand{\Ker}{\mathrm{Ker}}
 \newcommand{\pr}{\mathrm{pr}}
 \newcommand{\Spec}{\mathrm{Spec}}

 \newcommand{\Gr}{\mathrm{Gr}}

\renewcommand{\M}{\mathcal{M}}
\newcommand{\W}{\mathbf{W}}

\newcommand{\frob}{\mathrm{frob}}

\newcommand{\EExt}{\mathbf{Ext}}

\renewcommand{\Sym}{\mathrm{Sym}}

\newcommand{\fonction}[5]{\begin{array}{cccc}
#1: & #2 & \longrightarrow & #3 \\
 & #4 & \longmapsto & #5 \end{array}}
 
\newcommand{\fonctionnoname}[4]{\begin{array}{ccc}
 #1 & \longrightarrow & #2 \\
 #3 & \longmapsto & #4 \end{array}}

\newtheorem*{theorem*}{Theorem}
 \theoremstyle{plain}
 \newtheorem{thm}{Theorem}[section]
 \newtheorem{defi}[thm]{Definition}
 \newtheorem{prop}[thm]{Proposition}
 \newtheorem{lem}[thm]{Lemma}
 \newtheorem{coro}[thm]{Corollary}

 \theoremstyle{remark}
 \newtheorem{rem}[thm]{Remark}
 \newtheorem{qu}[thm]{Question}

 \newtheorem{ex}[thm]{Example}

\newtheorem*{mot}{Motto}

    \newtheorem{nota}[thm]{Notation}

 \newenvironment{dem}{{\bf Proof.}}{\hfill$\square$}

\date{}

\newif\ifquoteopen
\catcode`\"=\active 
\DeclareRobustCommand*{"}{%
   \ifquoteopen
     \quoteopenfalse ''%
   \else
     \quoteopentrue ``%
   \fi
}
\begin{document}
\quad
\vspace{1cm}
\begin{center}
{\Large \textsc{Lifting vector bundles to Witt vector bundles}}
\end{center}
\vspace{1cm}
{\normalsize \textsc{\noindent Charles De Clercq\!\footnote{Partially supported by French ministries of Foreign Affairs and of Education and Research (PHC Sakura- New Directions in Arakelov Geometry).}, Mathieu Florence and Giancarlo Lucchini-Arteche\!\footnote{Partially supported by Conicyt/ANID via Fondecyt Grants 11170016 and 1210010 and PAI Grant 79170034.}}}
\vspace{0.3cm}




\address{Charles De Clercq, Equipe Topologie Alg\'ebrique, Laboratoire Analyse, G\'eom\'etrie et Applications, Universit\'e Sorbonne Paris Nord, 93430 Villetaneuse, France.}
\address{Mathieu Florence, Universit\'e Paris Cit\'e and Sorbonne Universit\'e, CNRS, IMJ-PRG, F-75005 Paris, France }
\address{Giancarlo Lucchini Arteche, Departamento de Matem\'aticas, Facultad de Ciencias, Universidad de Chile, Las Palmeras 3425, \~Nu\~noa, Santiago, Chile.}

\Abstract
Let $X$ be a scheme. Let $r \geq 2$ be an integer. Denote by $\W_r(X)$ the scheme of Witt vectors of length $r$, built out of $X$. We are concerned with the question of extending  (=lifting) vector bundles on $X$, to vector bundles on $\W_r(X)$---promoting a systematic use of Witt modules and Witt vector bundles.
To begin with, we investigate two elementary but significant cases, in which the answer to this question is positive: line bundles, and the tautological vector bundle of a projective bundle over an affine base. We then  offer a simple (re)formulation of classical results in deformation theory of smooth varieties over a field $k$ of characteristic $p>0$, and extend them to \emph{reduced} $k$-schemes. Some of these results were recently recovered, in another form, by Stefan Schr\"oer. As an application, we prove that the tautological vector bundle of the Grassmannian $\Gr_{\F_p}(m,n)$ does not extend to  $\W_2(\Gr_{\F_p}(m,n))$, if $2 \leq m \leq n-2$. To conclude, we establish a connection to the  work of Zdanowicz, on non-liftability of some projective bundles. 

\tableofcontents

\section{Introduction}

Let $p$ be a prime, and let $X$ be a scheme, not necessarily of characteristic $p$. For any $r \geq 1$, denote by $\W_r(X)$ the scheme of $p$-typical Witt vectors of length $r$, built out of $X$. Let $V$ be a vector bundle over $X$. This paper mainly deals with

\begin{qu}\label{ques main}
Is $V$ the restriction to $X$ of a vector bundle defined over $\W_r(X)$?
\end{qu}
To fix ideas, until the end of the introduction, assume that $X$ is an $\F_p$-scheme.
The closed immersion $X \hookrightarrow \W_r(X)$ can then be thought of as a universal thickening of $X$, of characteristic $p^r$. Extending $V$ to $\W_r(X)$ is, in a sense made precise in Section \ref{sec TautoLift}, the "universal" deformation problem for $V$.\\

In order to tackle Question \ref{ques main}, we introduce in Section \ref{sec Witt modules} the notions of \emph{Witt modules} and \emph{Witt vector bundles}. Focusing on these objects (which are not at all new) provides a novel viewpoint on some problems in deformation theory. Rather than ``extending $V$ to a vector bundle over $\W_r(X)$'', we  say ``lifting $V$ to a $\W_r$-bundle \emph{over $X$}''---both formulations being of course equivalent. Doing so, the base scheme $X$ remains unchanged, while considering a larger class of sheaves of modules over it.  Section \ref{sec Witt modules} is devoted to laying foundations for this theory.\\

Section \ref{sec TeichLift} elaborates on a positive elementary answer to Question \ref{ques main}, for line bundles. The starting point here, is that every line bundle $L$ admits a natural lift to a $\W_r$-bundle: its $r$-th Teichm\"uller lift $\W_r(L)$. Section \ref{sec TautoLift} deals with lifting  the tautological vector bundle on $Y=\P_X(V)$, the projective space of a vector bundle $V$ over an affine $\F_p$-base $X$,  respecting the  tautological exact sequence. The answer is positive again---see Theorem \ref{PropTautoLift}.\\

In Section \ref{sec deformation}, we use the point of view of Witt modules to provide a unified functorial formulation of three questions in deformation theory:
\begin{enumerate}
    \item Let $X$ be a scheme over a perfect field $k$, of characteristic $p$. Lift $X$ to a flat $\W_2(k)$-scheme $X_2$ (a problem made famous by the article \cite{DI}).\vspace{0.2cm}
    \item Same as (1), with the strong extra constraint, that the Frobenius of $X$ lifts to $X_2$ (see the nice paper \cite{MS}).\vspace{0.2cm}
    \item Lift a vector bundle over $X$ to a vector bundle over $\W_2(X)$---the main topic of this paper.
\end{enumerate}
To achieve this, we describe equivalences of categories interpreting these questions in the common framework of Witt modules, see Propositions \ref{EquLiftW2}, \ref{EquLiftFrobW2} and \ref{EquW2}. These hold for any \textit{reduced} $X$, thus providing a gain of generality with respect to classical results, which require that $X$  be \textit{smooth} over $k$.

Several applications (some classical, some new) of these equivalences are given.
\begin{itemize}
    \item  Proof that (1) has an affirmative and functorial answer, for any Frobenius-split reduced $k$-scheme---see Corollary \ref{LiftFrobSplit} (for an introduction to Frobenius splitting, see \cite{MR}).
    \item Explicit $2$-extensions realizing the classical cohomological obstructions in deformation theory---see Corollaries \ref{explicitI} and \ref{explicitII}.
    \item In Theorem \ref{nonlift}, it is shown that the tautological vector bundle on the Grassmannian $\Gr(m,n)$ never lifts to a $\W_2$-bundle, if $2 \leq m \leq n-2$. Two different proofs are given.
    \item Finally, using the previous item, Section \ref{Appli2} provides a new proof, and a slight strenghtening,  of a result of Zdanowicz (\cite{Z}, Theorem 6.5).
\end{itemize}

\textbf{Acknowledgments.} We are grateful to  Bhargav Bhatt, Laurent Fargues, Ofer Gabber, Luc Illusie,  Boris Kunyavski\u{\i}, Marco Maculan, Matthew Morrow, Alexander Petrov, Stefan Schr\"oer, and anonymous referees, for valuable suggestions and enjoyable interactions.

\section{Conventions}
All schemes are assumed to be quasi-compact and quasi-separated. 

\section{Witt modules and \texorpdfstring{$\W_r$}{Lg}-bundles}\label{sec Witt modules}

 For each integer $r \geq 1,$ and for each commutative ring $A$, denote by $\W_r(A)$ the ring  of $p$-typical Witt vectors of $A$, of length $r$. Note that $\W_r$  itself might be considered as a ring scheme defined over $\Z$, isomorphic, as a scheme, to the $r$-dimensional affine space $\mathbb A^r$. For details on the construction of Witt vectors (adopting different viewpoints) we refer to \cite{BO}, \cite{Ha}, \cite{Le}, and \cite{Se1}.\\
 Fix a scheme $X$, not necessarily of characteristic $p$. Cover it by affine open subschemes $U_i=\Spec(A_i)$. Then, the affine schemes $\Spec(\W_r(A_i))$ glue, giving rise to $\W_r(X)$, the scheme of $p$-typical Witt vectors of length $r$, built out of $X$.

\subsection{Basic definitions}

The following definition is classical, see \cite{Se2}.

\begin{defi}
Let $r \geq 1$ be an integer. The association
\[U \mapsto \W_r(\mathcal O_X(U)),\]
defines a sheaf of commutative rings on $X$. Denote it by $\W_r(\mathcal O_X)$.
\end{defi}

As a sheaf of sets, $\W_r(\mathcal{O}_X)$ is represented by the $r$-dimensional affine space $\mathbb A^r_X$; it is thus also a sheaf for the fppf topology. If $1 \leq s < r$ are integers, denote by $$\pi_{r,s}: \W_r(\mathcal O_X) \to \W_s(\mathcal O_X)$$ the quotient, fitting into an exact sequence of sheaves of commutative rings on $X$ \[0 \to \W_{r-s}(\mathcal O_X) \xrightarrow{i_{r-s,r}} \W_r(\mathcal O_X) \xrightarrow{\pi_{r,s}} \W_s(\mathcal O_X) \to 0.\]

\begin{defi}
Let $r \geq 1$ be an integer. The association
\[U \mapsto \W_r(\mathcal O_X(U))^\times,\]
defines a sheaf of Abelian groups on $X$. It is given  by ($U$-points of) an affine and smooth $\Z$-group scheme,  denoted by $ \W_r^\times$.\\
\end{defi}

Note that  $\W_1 ^\times= \G_m$. The Teichm\"uller representative is given by the formula
\[\fonction{\tau_r}{\W_{1} ^\times}{\W_{r} ^\times,}{x}{(x,0,0, \ldots,0).}\]
We denote this homomorphism of $\Z$-groups schemes by $\tau$, if the dependence in $r$ is clear. It is a natural splitting of the quotient map $\W_{r} ^\times \to \G_m. $\\

\begin{rem}
Let $R$ be an $\F_p$-algebra. For all $t \in R$, the following formula holds:
\[ p \tau_{r+1}(t) = i_{r,r+1} (\tau_r(t^p)) \in \W_{r+1}(R).\]
This can be used to check that, for all $x,y \in \W_r(R)$,
\[i_{r,r+1}(x)i_{r,r+1}(y)=p\, i_{r,r+1}(xy),\]
where multiplication on the left (resp. right), is that of $\W_{r+1}(R)$ (resp. $\W_r(R)$).\\
\end{rem}

Let $r \geq 1$ be an integer. There is an isomorphism of affine $\Z$-group schemes
\begin{align*}
    \G_m \times (1+\W_r)^\times &\overset{\sim}{\longrightarrow}\W_{r+1} ^\times\\
    (a,1+x) & \longmapsto \tau_{r+1}(a)(1+x).
\end{align*}
(Strictly speaking, one should write $(1+i_{r,r+1}(\W_r))^\times $ instead of $(1+\W_r)^\times$.) This formula, as well as the next ones, are given at the level of functors of points.

As a morphism of linear algebraic $\F_p$-groups, the logarithm is well-defined: 
\begin{align*}
\log: (1+\W_r)^\times  &\to (\W_r,+),\\
(1-x) &\longmapsto -x- \frac p 2 x^2 -  \frac {p^2} 3 x^3 - \frac {p^3} 4 x^4 - \ldots
\end{align*}
If $p$ is odd, it is an isomorphism of linear algebraic $\F_p$-groups, with inverse
\begin{align*}
\exp: (\W_r,+) &\to (1+\W_r)^\times,\\
t &\longmapsto 1+ t+ \frac p {2!} t^2 +  \frac {p^2} {3!} t^3 + \ldots \end{align*}
The sums occuring above are well-defined, since $\frac {p^i} i $ (resp. $\frac {p^i} {i!} $ when $p$ is odd) is a $p$-adic integer, for all $i \geq 0$. Moreover,  these sums are finite, because $x$ is nilpotent.

If $p=2$, the logarithm is an isogeny of degree two, with kernel $\{ 1,-1\}$. Over $\F_2$, the algebraic group $(1+\W_r)^\times$ is then isomorphic to the middle term of the pull-back of the exact sequence
\[ 0 \to \W_{r-1} \to   \W_{r} \to \G_a \to 0, \]
by the Lang isogeny
$$(\frob-\Id): \G_a  \to \G_a.$$
Thus, from the point of view of linear algebraic groups, the multiplicative group scheme $\W_{r+1}^\times$, over $\F_p$, is nothing new. \\However, over $\Z$, the morphisms $\log$ and $\exp$ above are not defined, and $\W_{r+1}^\times$ is much more intriguing, as the following remark shows.

\begin{rem}\label{ToreToWitt}
Put
\[\G_{\mathrm{a}/\mathrm{m}}:=\Ker(\W_2^\times \to \W_1^\times).\]
it is a smooth affine group scheme over $\Z$. Its generic fiber $\G_{\mathrm{a}/\mathrm{m}} \times_\Z \Q$ is isomorphic to $\Gm$, whereas its special fiber $\G_{\mathrm{a}/\mathrm{m}} \times_\Z \F_p$ is isomorphic to $\Ga$. More generally, consider
$$D_r:=\Ker (\W_{r+2} ^\times \to \W_2 ^\times);$$   Its generic fiber is a split algebraic torus of dimension $r$, whereas its special fiber is isomorphic to $(\W_r,+)$ (also for $p=2$). These are simple examples of deformations of (additive groups of) Witt vectors to (split) algebraic tori, in the spirit of \cite{TO}.\\
Let $\Lambda$ be a discrete valuation ring, of mixed characteristic $(0,p)$ (e.g. $\Lambda=\Z_p$). Over $\Lambda$, one may attempt to twist $D_r$ (by a suitable $\Aut(D_r)$-torsor), in order to obtain more interesting deformations of Witt vectors, to non-split algebraic tori. We will not follow this track here, but the interested reader is welcome to do so.
\end{rem}

We now  introduce Witt modules and $\W_n$-bundles, and  their first properties.

\begin{defi}\label{WtFaffine}
Assume that $X=\Spec(A)$ is affine. Let $r \geq 1$ be a positive integer. Let $M$ be a $\W_r(A)$-module. The formula $$ U \mapsto M \otimes_{\W_r(A)} \W_r(\mathcal O_X(U))$$ defines a presheaf (for the Zariski topology) on $X$. Denote by $\tilde M$ the associated sheaf. It is a sheaf of $\W_r(\mathcal O_X)$-modules.
\end{defi}

\begin{defi}[Witt module]
Let $r \geq 1$ be a positive integer. A Witt module of height $r$ over $X$ is a sheaf of $\W_r(\mathcal O_X)$-modules, which is locally isomorphic to a sheaf of the shape $\tilde M$ (cf. Definition \ref{WtFaffine}).
\end{defi}

\begin{defi}[$\W_r$-bundle]
Let $r,n \geq 1$ be two positive integers. A $\W_r$-bundle over $X$ of rank $n$ is a Witt module of height $r$, locally isomorphic to $\W_r(\mathcal O_X)^n$.
\end{defi}

\begin{rem}
A Witt module of height $1$ over $X$ is a quasi-coherent $\mathcal O_X$-module, while a a $\W_1$-bundle is a vector bundle over $X$. Witt modules of height $r$ (resp. a $\W_r$-bundle) are simply modules (resp. vector bundles) over the scheme $\W_r(X)$.
\end{rem}

\begin{nota}
Let $V_r/X$ be a $\W_r$-bundle over $X$. Denote by $ V_r^\vee$ the sheaf
\[\underline{\Hom}_{\W_r(\mathcal O_X)-\mathrm{Mod}} ( V, \W_r(\mathcal O_X));\]
it is a $\W_r$-bundle over $X$.

If $s \leq r$ is a positive integer,  define the $\W_s$-bundle 
\[V_r\otimes_{\W_r} \W_s\]
as the sheaf associated to the presheaf 
\[U \mapsto  V (U) \otimes_{\W_r(\mathcal O_X(U))} \W_s(\mathcal O_X(U)).\] 
\end{nota}

\subsection{Frobenius twists, for $\F_p$-schemes}\label{SectFrob}
Keep notation of the previous paragraph. Assume moreover, that $X$ is an $\F_p$-scheme. Denote by
\[ \frob_X: X \to X,\]
the (absolute) Frobenius of $X$, given by raising functions to their $p$-th power. By functoriality, it induces an arrow
\[\W_r(X) \to \W_r(X),\]
for each $r \geq 1$, which we still denote by $\frob_X$, or just by $\frob$.

If $M$ is a $\W_r(\mathcal O_X)$-module, we set, for each $n \geq 1$, \[M^{(n)}:=(\frob^n)^*(M).\]
Considering $M$ as a $\W_s(\mathcal O_X)$-module, for some $s >r$, does not affect the formation of $M^{(n)}$. This follows from the commutative diagram
\[ \xymatrix{ \W_r(X) \ar[r]^{\frob^n} \ar[d] & \W_r(X)  \ar[d]\\ \W_s(X) \ar[r]^{\frob^n} & \W_s(X), }\]
where the vertical arrows are the natural immersions. The notation $M^{(n)}$ thus makes sense for Witt modules.

For $s >r$, let  $V_s$ be a $\W_s$-bundle. Then  $V_s^{(n)}$ is a $\W_s$-bundle as well. Moreover, there is an exact sequence of Witt modules on $X$
\[\mathcal E_r(V_s): 0 \to (\frob^r)_*  (V_{s-r}^{(r)}) \to V_s \stackrel {\pi_{s,r}} \to V_{r} \to 0. \]
If $s=r+1$, and if $X$ is regular, then $(\frob^r)_*  (V_{s-r}^{(r)})$ is a vector bundle on $X$.

Recall that, in the case where $r=1$ and $L$ is a line bundle, $L^{(1)}$ is naturally isomorphic to $L^{\otimes p}$. This fundamental fact does not extend to vector bundles of dimension greater than $2$, nor to  $\W_2$-line bundles.

\subsection{Frobenius twists, general case.}\label{SectFrob2}
Consider the morphism of schemes
\[\Frob: \W_r(X) \to \W_{r+1}(X),\]
that functorially arises from the morphism of ring schemes over $\Z$ (cf. \cite[Definition 1.1]{DK})
\[\Frob: \W_{r+1} \to \W_r.\]
Note that, for $r=1$, it is given by the Witt polynomial
\[(x_0,x_1) \mapsto x_0^p +px_1.\]
For a $\W_{r+1}$-module $M_{r+1}$ over $X$, $\Frob^*(M_{r+1})$ is a $\W_r$-module over $X$. Note in particular the length shift by $(-1)$ in this general notion of Frobenius pull-back.\\
If $X$ is an $\F_p$-scheme, then $\Frob$ equals 
the composite
\[\W_r(X) \xrightarrow{\frob } \W_r(X)  \stackrel {nat} \hookrightarrow  \W_{r+1}(X),\]
where $nat$ is the natural nilpotent immersion. Thus, 
$\Frob^*(M_{r+1})$ depends only on  $M_r:=M_{r+1} \otimes_{\W_{r+1}} \W_r$, and  coincides with $M_r^{(1)}=\frob^*(M_r)$.\\
To avoid confusion, we reserve the notation $(\cdot)^{(1)}$ for the Frobenius pull-back of a $\W_r$-bundle over an $\F_p$-scheme $X$, via $\frob: \W_r(X) \to \W_r(X)$.\\
For $s >r$ and $V_s$ a $\W_s$-bundle, $(\Frob^r)^* (V_s)$ is a $\W_{s-r}$-bundle and there is an exact sequence of Witt modules on $X$
\[\mathcal E_r(V_s): 0 \to (\Frob^r)_* (\Frob^r)^* (V_s) \to V_s \stackrel {\pi_{s,r}} \to V_{r} \to 0. \] 

\subsection{Yoga of extensions}
In this and all subsequent sections, we will make constant use of the notions of extensions of Witt modules, endowed with the following classical notions: Baer sum, pullback, pushforward, change of the base. We briefly recall these, without proofs.

\begin{defi}(Extensions and operations on them)\\
Let $\mathcal A$ be an Abelian category. Let
\[ \mathcal E: 0 \to A \to B \to C \to 0,\]
be an exact sequence in $\mathcal A$, thought of as an extension of $C$ by $A$.

Let $f:A \to A'$ be a morphism in $\mathcal A$. Then $f_*(\mathcal E)$ is the unique extension fitting in a commutative diagram
\[\xymatrix{\mathcal E: 0 \ar[r] & A \ar[r] \ar[d]^f & B \ar[r] \ar[d] &C \ar[r] \ar@{=}[d] & 0 \\ f_*(\mathcal E):  0 \ar[r] & A' \ar[r]  & D \ar[r]  &C \ar[r] & 0.}\]
We refer to $f_*(\mathcal E)$ as the pushforward of $\mathcal E$ by $f$.

Let $g:C' \to C$ be a morphism in $\mathcal A$. Then $g^*(\mathcal E)$ is the unique extension fitting in a commutative diagram
\[\xymatrix{g^*(\mathcal E): 0 \ar[r] & A \ar[r] \ar@{=}[d] & E \ar[r] \ar[d] &C' \ar[r] \ar[d]^g & 0 \\ \mathcal E:  0 \ar[r] & A \ar[r]  & B \ar[r]  &C \ar[r] & 0.}\]
We refer to $g^*(\mathcal E)$ as the pullback of $\mathcal E$ by $g$.
\end{defi}

Let
\[ \mathcal E': 0 \to A' \to B' \to C'\to 0,\]
be another extension in $\mathcal A$. Then one can easily show that the data of a commutative diagram  \[\xymatrix{\mathcal E: 0 \ar[r] & A \ar[r] \ar[d]^f & B \ar[r] \ar[d] &C \ar[r] \ar[d]^g & 0 \\ \mathcal E':  0 \ar[r] & A' \ar[r]  & B' \ar[r]  & C' \ar[r] & 0,}\]
is equivalent to an isomorphism (of extensions of $C$ by $A'$) $f_*(\mathcal E) \stackrel \sim \to g^*(\mathcal E')$.

\begin{rem}
Let
\[\mathcal E:  0 \to \mathcal O_X \to B \to C \to 0,\]
be an extension of (quasi-coherent) modules over a scheme $X$, of characteristic $p$. In the sequel, the notation $\frob_*(\mathcal E)$ may refer either to the pushforward
\[\frob_*(\mathcal E):  0 \to \frob_*(\mathcal O_X) \to \tilde{B} \to C \to 0,\]
in the sense of the definition above, or to the extension  \[\frob_*(\mathcal E):  0 \to \mathcal \frob_*(\mathcal O_X) \to \frob_*(B) \to \frob_*(C) \to 0,\]
obtained by applying the exact functor (on $\mathcal O_X$-modules) $\frob_*$. This ambiguity is nothing serious: the meaning of $\frob_*(\mathcal E)$ is always clear from the context. The same goes for the notation $\frob^*(\mathcal E)$.
\end{rem}

\subsection{The arrow $\kappa$.}\label{sec kappa}
Assume given an extension of Abelian groups \[ 0 \to A \to B \stackrel \pi \to C \to 0, \] where $A$ and $C$ are killed by $p$ (hence $B$ is killed by $p^2)$. By assumption, for any $c\in C$, setting $$\kappa(c):= pb,$$ where $b$ is such that $\pi(b)=c$ gives rise to a well-defined morphism $$\kappa:C\to A.$$ Furthermore, $\kappa$ is trivial if and only if $pB=0$, and it is an isomorphism if and only if $B$ is a lift of the $\F_p$-vector space $C$ to a free $(\Z /p^2 \Z)$-module.

One can clearly ``sheafify'' the construction of the arrow $\kappa$ in the following way.\\Let $\mathcal S$ be a site, and let
\[ \mathcal E: 0 \to A \to B \stackrel \pi \to C \to 0,\]
be an exact sequence of sheaves of Abelian groups on $\mathcal S$. Assume given an endomorphism $f\in \Hom_{\mathcal S}(  B , B)$, leaving $A$ stable. It thus induces an endomorphism of $C$, still denoted by $f$. Assuming that $f(A)=f(C)=0$, define a homomorphism
\[ \kappa=\kappa_{\mathcal{E},f} \in \Hom_{\mathcal S}(  C , A),\]
through the following process. Let $Y$ be an object of $\mathcal S$. For any section $s \in C(Y)$, pick a covering $(Y_i)$ of $Y$ such that, for each $i$,  the restriction $s_i=s_{\vert Y_i}$ lifts to $t_i \in B(Y_i)$. Then $u_i:=f(t_i)$ belongs to $A(Y_i)$. The $u_i$ hence glue to give an element in $ A(Y)$, depending only on $s$. Denote it by $\kappa(s).$ Note that $\kappa=0$ if and only if $f=0$.

Now if $\mathcal S$ is either the Zariski, \'etale or fppf  site, big or small,  of a scheme $X$ of characteristic $p$, consider the exact sequence of Witt modules over $X$
\[ \mathcal E: 0 \to \Frob_*(\mathcal O_X) \to \W_{2}(\mathcal O_X ) \stackrel \pi \to  \mathcal O_X \to 0.\]
and take $f$ to be multiplication by $p$. Then, the homomorphism $\kappa_{\mathcal E,p}$ is given by the Frobenius
\[\Frob: \mathcal O_X \stackrel {x \mapsto x^p} \to  \Frob_*(\mathcal O_X).\]
Now assume that $X$ is defined over a perfect field $k$ and that $X/k$ admits a lift $X_2$, flat over $\W_2(k)$. Recall that a  $\W_2(k)$-module is flat, if and only if it is free. Consider the exact sequence
\[ 0 \to p\mathcal O_{X_2}  \to \mathcal O_{X_2}  \stackrel {can} \to  \mathcal O_X  \to 0,\] whose kernel is a square-zero ideal of $ \mathcal O_{X_2} $. If $f$ is multiplication by $p$, then $\kappa$ is an isomorphism.

\subsection{Lifting Witt vector bundles}
In this section, $X$ is an arbitrary scheme.
\begin{defi}
Let $ r \geq 1$ be  an integer. Let $V_r/X$ be a $\W_r$-bundle.
\begin{itemize}
    \item Let $s >r$ be an integer. A lifting of $V_r$ to a $\W_s$-bundle is the data of a $\W_s$-bundle $ V_s$ on $X$, together with an isomorphism (of $\W_r$-bundles) $$ f_r:  V_{s}\otimes_{\W_s} \W_r \stackrel \sim \to  V_{r}.$$ If $X$ is an $\F_p$-scheme, one also says that $V_s$ is a lift of $V_r$ to $p^s$-torsion.\\
    
    \item { A complete lifting of $V_r$ is the data, for each $s> r$, of a $\W_s$-bundle $ V_s$ on $X$, together with isomorphisms (of $\W_s$-bundles) $$ f_s:  V_{s+1}\otimes_{\W_{s+1}} \W_s  \stackrel \sim \to  V_{s}.$$}
\end{itemize}
\end{defi}

The notion of a complete lifting of $V_r$ is the same as that of a lift of $V_r$ to a $\W_\infty$-bundle, where $\W_\infty =\varprojlim \W_r$. In the case where $r=1$ and $X= \Spec(\F_p)$, a complete lifting for $V_1$ is then the data of a free $\Z_p$-module of finite rank $ V_\infty,$ together with an isomorphism of $\F_p$-vector spaces $V_\infty/p \stackrel \sim \to V_1$.
\begin{rem}
The existence of (unrelated) lifts of $V_r$ to $p^s$-torsion, for each $s >r$, need  not  imply the existence of a complete lifting of $V_r$.     
\end{rem}

\begin{rem}
    
Assume that $X$ is an $\F_p$-scheme.
Liftability of a $\W_r$-bundle $V_r$ to a $\W_{r+1}$-bundle is then equivalent to the vanishing of a natural class in the coherent cohomology group $H^2(X,\End(V_1)^{(r)})$. 
More generally, assume that $r \leq s \leq 2r$. Then, lifting $V_r$ to $p^s$-torsion is obstructed by a class in  $H^2(X,\End(V_{s-r})^{(r)})$, where $V_{s-r}:=V_r \otimes_{\W_r} \W_{s-r}$. One checks this using the exact sequence $\mathcal E_r(V_s)$ of Section \ref{SectFrob}.
\end{rem}

\begin{rem}
    In general,   liftability of a vector bundle on $X$ to a $\W_2$-bundle, is a highly non-abelian problem.
\end{rem}
\section{The Teichm\"uller representative lifts line bundles}\label{sec TeichLift}
The Teichm\"uller representative  yields natural lifts for line bundles. This construction is detailed below and works over an arbitrary scheme $X$. However, some results are specific to $\F_p$-schemes---the case of interest for applications later on.

\subsection{Lifts of line bundles.}
\begin{defi}
Let $ L$ be a line bundle over an arbitrary scheme $X$. Denote by \[P:=\underline{\mathrm{Isom}}(L,\mathcal O_X)=\mathbb A(L)-\{0\},\]
the corresponding $\G_m$-torsor over $X$. Let $r \geq 1$ be an integer. Then
\[P_{r}:=(\tau_{r})_*(P),\]
is a $\W_{r}^\times$-torsor over $X$.

Twisting the trivial invertible $\W_{r}$-module $\W_{r}(\mathcal O_X)$ by $P_{r}$ yields an invertible $\W_{r}$-module. We denote it by $\W_{r}(L)$. It is the $r$-th Teichm\"uller lift of $L$.
\end{defi}

\begin{rem}
Note that, for a group scheme $G$, twisting $G$-schemes by $G$-torsors usually requires extra assumptions---e.g. quasiprojectivity. Here, for $G=\W_r ^\times$, no assumption is needed since line bundles are locally trivial for the Zariski topology.
\end{rem}

\begin{rem}
There is a natural isomorphism of $\W_r$-line bundles \[\Frob^*(\W_{r+1}(L))=\W_r(L^{\otimes p}),\]
arising functorialty from the fact that the composite \[\W_1^\times \xrightarrow{\tau_{r+1}} \W_{r+1}^\times  \xrightarrow{\Frob} \W_r^\times,\]
coincides with
\[\W_1^\times \xrightarrow{x \mapsto x^p} \W_1^\times  \xrightarrow{\tau_r} \W_r^\times.\] 
\end{rem}

The sequence $(\W_r( L))_{r \geq 1}$, together with the natural isomorphisms $$\W_{r+1}(L) \otimes_{\W_{r+1}} \W_r \simeq \W_r( L),$$ form a natural complete lifting for $L$, which is functorial in the following sense.

Pick integers  $s > r \geq 1$. Consider the exact sequence of Witt modules on $X$
\[0 \to  (\Frob^r)_*(\W_{s-r} (\mathcal  O_X)) \xrightarrow{i=i_{s-r,s}} \W_{s} (\mathcal O_X) \stackrel {\pi=\pi_{s,r}} \to  \W_{r} (\mathcal O_X) \to 0.\]
Written as such, the arrows $i$ and $\pi$ are $\G_m$-equivariant, where  $\G_m$ acts on  $ \W_{r}(\mathcal O_X)$ via the multiplicative section $\tau_{r}$.

Twisting by the $\G_m$-torsor $P$  associated to a line bundle $L$ over $X$, one gets an exact sequence of Witt modules on $X$
\[0 \to (\Frob)^r_*(\W_{s-r}(L^{\otimes p^r} )) \to  \W_{s} (L) \xrightarrow{\pi_L=\pi_{ L,s,r}} \W_{r} (L) \to 0,\]
which coincides with the sequence $\mathcal E_r (\W_{s}(L))$ from the end of Section \ref{SectFrob2}.\\

We describe now the space of lifts of a line bundle to a $\W_2$-line bundle, over an $\F_p$-base.

\begin{prop}\label{LiftToW2}
Assume that $X$ is an $\F_p$-scheme.
Let $X$ be a scheme and $L_1$ be a line bundle over $X$. Lifts of $L_1$ to a $\W_2$-bundle $L_2$ over $X$ are in natural bijection with $\G_a$-torsors over $X$. Hence, the set of isomorphism classes of such $L_2$'s  is in natural bijection with $H^1(X, \mathcal O_X)$.
\end{prop}

\begin{dem}
Let $L_2$ be a lift of $L_1$, to a $\W_2$-bundle over $X$. It fits into the extension of $\W_2$-modules over $X$
\[\mathcal E_1(L_2): 0 \to   \frob_*(L_1^{\otimes p})=\frob_*(L_1^{(1)}) \to L_2 \to L_1 \to 0.\]
In particular, taking $L_2=\W_2(L_1)$ yields
\[\mathcal E_1(\W_2(L_1)): 0 \to   \frob_*(L_1^{\otimes p}) \to \W_2(L_1) \to L_1 \to 0.\]
The $\kappa$ arrow (Section \ref{sec kappa}) of both extensions is $\frob_{L_1}$. Form the Baer difference 
\[\Delta(L_2):=\mathcal E_1(L_2)-\mathcal E_1(\W_2(L_1)): 0 \to \frob_*(L_1^{\otimes p})  \to D(L_2) \to L_1 \to 0.\]
Since the arrow $\kappa$ commutes with Baer sum, it is trivial for $\Delta(L_2)$. Its middle term $D(L_2)$ is thus a quasi-coherent $\mathcal O_X$-module, so that $\Delta(L_2)$ is actually an extension of quasi-coherent $\mathcal O_X$-modules. Using the adjunction between $\frob_*$ and $\frob^*$, $\Delta(L_2)$ gives rise to an extension of vector bundles over $X$
\[\epsilon(L_2): 0 \to L_1^{\otimes p} \to E(L_2) \to L_1^{\otimes p} \to 0.\]
Applying $\cdot \otimes L_1^{\otimes -p}$ yields an extension of vector bundles over $X$
\[\alpha(L_2):=\epsilon(L_2) \otimes  L_1^{\otimes -p}: 0 \to \mathcal O_X \to A(L_2) \to \mathcal O_X \to 0;\]
that is to say, a $\G_a$-torsor over $X$.
 
This construction is reversible, providing an equivalence between liftings of $L_1$ to a $\W_2$-line bundle $L_2$, and  $\G_a$-torsors over $X$, which are classified by $H^1(X, \mathcal O_X)$.
\end{dem}

\begin{rem}
Alternatively, one can prove Proposition \ref{LiftToW2} in the following way. Recall that $H^1(S, \W_r^\times)$ classifies $\W_r$-line bundles over $S$. Consider the natural isomorphism of $\F_p$-group schemes  \[ \G_a  \times \G_m \stackrel \sim \to \W_2 ^\times , \] and  apply $H^1(S,.)$, to get an isomorphism \[ H^1(S,\mathcal O_S) \oplus  H^1(S,\G_m) \stackrel \sim \to H^1(S,\W_2 ^\times). \]
This point of view generalizes to arbitrary schemes as follows. Over $\Z$, consider the affine group scheme
$\G_{a/m}$ defined in Remark \ref{ToreToWitt}
Let $L_1$ be a line bundle over an arbitrary scheme $X$. Then one can prove that isomorphism classes of lifts of $L_1$ to a $\W_2$-line bundle $L_2$ are in natural bijection with $H_{\mathrm{Zar}}^1(X,\G_{a/m})$.
\end{rem}

\subsection{The Teichm\"uller section for line bundles}\label{TeichWitt}
Keep notation of the preceding paragraph.
The surjection $\pi: \W_{s} \to \W_r$ for $s>r$ has a natural multiplicative scheme-theoretic section,
\begin{align*}
    \tau(=\tau_{r,s}): \W_r  &\to \W_{s},\\
    (x_1, \ldots, x_r) &\longmapsto (x_1, \ldots, x_r,0,\ldots,0).
\end{align*}
This section is $\G_m$-equivariant, for the $\G_m$-action on $ \W_r$ (resp. on  $ \W_{s}$) given by $\tau_{1,r}$ (resp. $\tau_{1,s}$). Twisting by the torsor $P$ associated to $L$ then yields a natural, scheme-theoretic section of $\pi_{ L}$, its Teichm\"uller section, denoted by
\[\tau_{L}: \W_r(L)  \to \W_{s}(L).\]
Note that $\tau_{L}$ is of course not additive. Taking global sections, one gets  exactness of the sequence of abelian groups
\[ 0 \to H^0(X,  (\Frob^r)^*(\W_s(L))) \to H^0(X,   \W_{s}(L) ) \to H^0(X,\W_{r}(L) ) \to 0.\]

Furthermore, assume that $G$ is a finite group  acting on $X$  by scheme automorphisms, and semi-linearly on $L$; in other words, that $L$ is a $G$-linearized line bundle. Then, the natural map on $G$-invariant sections
\[ H^0(G, H^0(X,\W_{s}(L) )) \to H^0(G,H^0(X,\W_{r}(L) ))\]
is onto as well. Indeed, the Teichm\"uller section $\tau_{ L}$ is  then $G$-equivariant.

\section{Lifting the tautological vector bundle of a projective space} \label{sec TautoLift}

\subsection{Two  facts about lifting vector bundles}
The following Lemma asserts that being liftable is invariant under twists.

\begin{lem}
Let $X$ be a scheme, let $V$ be a vector bundle on $X$, and let $L$ be a line bundle on $X$.  Let $r \geq 2$ be an integer. Then, $V$ lifts to a $\W_r$-bundle if and only if $V \otimes L$ does.
 \end{lem}

\begin{dem}
Let $V_r$ be a lift of $V=V_1$ to a $\W_r$-bundle over $X$. By Section \ref{sec TeichLift}, the line bundle $L$ admits the   lift $\W_r(L)$. Then, $V_r \otimes \W_r(L)$ is a lift of $V\otimes L$.
\end{dem}

In some sense,  liftability of vector bundles to Witt vector bundles is ``the hardest'' existence problem in the deformation  theory of vector bundles. In characteristic $p$, this is made precise in Lemma \ref{LemLiftUniv} below.

\begin{defi}
Let $X$ be a scheme and let $i: X \hookrightarrow Y$ be a closed immersion.

We say that $i$ is a $p$-elementary thickening of $X$, if the sheaf of ideals $I \subset \mathcal O_Y$ defining $X$ satisfies $I^p=pI=0$.

Let $r \geq 1$ be an integer. We say that $i$ is a $p$-thickening of depth $r$ of $X$ if it can be written as a composite
\[s:X=X_0  \hookrightarrow X_1 \hookrightarrow \ldots  \hookrightarrow X_{r-1} \hookrightarrow X_r=Y,\]
of $r$ $p$-elementary  thickenings.
\end{defi}

\begin{ex}
If $X$ is a smooth variety over a perfect field $k$ of characteristic $p$, and if  $X_{r+1}$ is a lift of $X$ to a scheme  $X_{r+1}$, flat over $\W_{r+1}(k)$, then the closed immersion $X \hookrightarrow X_{r+1}$ is  a $p$-thickening of depth $r$.
\end{ex}

\begin{lem}\label{LemPrepW}
   Let $X$ be an arbitrary scheme and let $i: X \hookrightarrow Y$ be a closed immersion. Then $i$ is a $p$-thickening of  depth $r$ of $X$, if and only if the sheaf of ideals $I \subset \mathcal O_Y$ defining $X$ satisfies $p^{r-i} I ^{p^i}=0,$ for every $i=0, \dots, r$.
\end{lem}

\begin{dem}
 Assume that $I$ satisfies the conditions of the Lemma. For  $j=0,\ldots,r$, define the sheaf of ideals $$ J_j:= (I^{p^j}+pI^{p^{j-1}}+p^2I^{p^{j-2}}+\ldots +p^j I) \subset \mathcal O_Y.$$  Let $Z_j \hookrightarrow Y$ be the closed subscheme defined by $J_j$. It is straightforward to check that $$(J_j^p+pJ_j) \subset J_{j+1} \subset J_j,$$  just using $p \geq 2$ (that $p$ is prime is not needed here). Hence, the closed immersion  $Z_j \hookrightarrow Z_{j+1}$ is a $p$-elementary thickening. To conclude, observe that $Z_0=Y$ and $Z_r=X$.  The converse implication is straightforward as well.  
\end{dem}

\begin{lem}\label{LemLiftUniv}
Assume that $X$ is an $\F_p$-scheme. Let $r \geq 1$ be an integer, and let $i: X \hookrightarrow Y$ be a $p$-thickening of $X$, of depth $r$. Let $V$ be a vector bundle on $X$. Assume that $V$ admits a lift to a $\W_{r+1}$-bundle on $X$. Then, $(\frob^r)^*(V)$ extends (via $i$) to a vector bundle on $Y$.
\end{lem}

\begin{dem}
Denote by
$$f: X \to \W_{r+1}(X),$$
the natural thickening. There is a natural morphism
\[ F_r: Y \to  \W_{r+1}(X),\]
such that
$$ F_r \circ i = f \circ \frob^r.$$
This fact is classical, but we could not find a suitable reference in the literature. Laurent Fargues kindly suggested a quick proof, which we now provide.

By glueing, it is straightforward to reduce to the affine case $X=\Spec(A)$. Let  $I \subset A$ be the  ideal defining $X$. By Lemma \ref{LemPrepW}, it satisfies, for every $i=0, \dots, r$,
\begin{equation}\label{eqn arg Fargues}
p^{r-i} I ^{p^i}=0.
\end{equation}
Consider the natural ring homomorphism given by the Witt polynomial
\begin{align*}
\Phi_r: \W_{r+1}(A) &\to 
A,\\
x:=(x_0,\ldots, x_r) &\longmapsto  x_0^{p^r}+ px_1^{p^{r-1}}+  p^2x_2^{p^{r-2}}+ \ldots  + p^rx_r.
\end{align*}
By \eqref{eqn arg Fargues}, $\Phi_r$ vanishes when all $x_i$ belong to $I$. Equivalently, it vanishes when $x$ belongs to the kernel of  the natural exact sequence of abelian groups
\[0 \to   \W_{r+1}(I) \to   \W_{r+1}(A) \xrightarrow{\pi}   \W_{r+1}(A/I) \to 0.\]
Thus, $\Phi_r$ factors via $\pi$, to a ring homomorphism $$\phi_r: \W_{r+1}(A/I) \to 
A,$$  proving the claim (set $F_r:=\Spec(\phi_r)$).

Via $F_r^*$, the existence of an extension of $V$ to (a vector bundle over) $\W_{r+1}(X)$ then implies that of an extension of $(\frob^{r})^*(V)$ to $Y$.
\end{dem}

From now on, and until the end of this section, 
 $X$ is a scheme of characteristic $p >0,$ and $V/X$ is  a vector bundle, of constant rank $n \geq 1.$ Denote by \[ f:\P(V)  \to X \] the associated projective bundle,  parametrizing \textit{quotient}  line bundles of $V$.

\begin{defi}
Denote by $\mathcal H_V$, or simply by $\mathcal H$ if the dependence in $V$ is clear, the tautological vector bundle on $\P(V)$.\\ There is the tautological exact sequence of vector bundles on $\P(V)$ \[\mathcal T (=\mathcal T_V):  0 \to \mathcal H \to f^*(V) \to \mathcal O(1) \to 0. \]
\end{defi}

\begin{thm}\label{PropTautoLift}
Assume that the $\F_p$-scheme $X$ is affine. Then $V$ lifts completely. Moreover, if we choose a complete lifting $(V_r)_{r \geq 1}$ of $V$, over $X$, then there exists a complete lifting  $(\mathcal H_r)_{r \geq 1}$ of $\mathcal H$, over $\P(V)$, such that each $\mathcal H_r$ fits into an exact sequence  of $\W_r$-bundles on $\P(V)$ \[\mathcal T_r :  0 \to \mathcal H_r \to f^*(V_r) \to \mathcal \W_r(\mathcal O(1)) \to 0, \] and such that these exact sequences are compatible. In other terms, for each $r \geq 1$, there is   a commutative diagram \[   \xymatrix{\mathcal T_{r+1} : 0\ar[r] & \mathcal H_{r+1} \ar[r] \ar[d] & f^*(V_{r+1} ) \ar[r]^-{\rho_{r+1}} \ar[d] &  \W_{r+1}(\mathcal O(1)) \ar[r] \ar[d]^{\pi_{r+1,r}} & 0  \\ \mathcal T_r : 0\ar[r]  & \mathcal H_{r} \ar[r]  & f^*(V_{r} ) \ar[r]^-{\rho_{r}}  &  \W_{r}(\mathcal O(1)) \ar[r] & 0,}\]
where the middle vertical arrow is the natural one.
\end{thm}
  
\begin{dem}
Assume that $V_r$ is a given lift of $V$, to a $\W_r$-bundle. The obstruction to lifting  $V_r$ to a $\W_{r+1}$-bundle $V_{r+1}$ lies in
\[\Ext^2_{\mathcal O_X-\mathrm{Mod}} ( V, (\frob^r)_* (V^{(r)}))=H^2(X,  \End(V)^{(r)}).\]
This cohomology group vanishes since $X$ is affine, whence the first claim.

Let us prove the second claim: assume that $\mathcal H_r$, together with the extension
\[\mathcal T_r :  0 \to \mathcal H_r \to f^*(V_r) \stackrel {\rho_r} \to \mathcal \W_r(\mathcal O(1)) \to 0, \]
has been constructed. There is a natural duality isomorphism
\[\underline{\Hom}_{\W_r(\mathcal O_X)-\mathrm{Mod}} (f^*(V_r),\W_r(\mathcal O(1)) ) = f^*(V_r) ^\vee(1) \stackrel \sim \to f^*(V_r ^\vee) (1). \]
The surjection $\rho_r$ thus corresponds to a global section \[ s_r \in H^0(\P(V),f^*(V_r ^\vee)(1)).\] One would like to lift it, through the epimorphism of the exact sequence (on $\P(V)$) $$ 0 \to   \frob^r_*(f^*(V^{(r)\vee})) (1) \to f^*(V_{r+1} ^\vee)(1)  \to f^*(V_r ^\vee)(1) \to 0.$$
The obstruction to do so is a class $$c \in  H^1(\P(V),  f^*(V^{\vee(r)}) (p^r))= 0. $$
(To get this vanishing result, use  $H^1_{\mathrm{Zar}}(\P(V),  \mathcal O (p^r)) =0$, together with the projection formula.)\\
Hence, $s_r$ can be lifted to a global section
\[ s_{r+1} \in H^0(\P(V),f^*(V_{r+1} ^\vee)(1)).\]
Dualizing, it corresponds to a homomorphism
\[\rho_{r+1}: f^*(V_{r+1} )  \to \W_{r+1}(\mathcal O(1) ), \] 
lifting $\rho_r$. Since $\rho_r$ is surjective, $\rho_{r+1}$ is surjective as well (using Nakayama's Lemma). Define the $\W_{r+1}$-bundle  $\mathcal H_{r+1} $ to be its kernel. Existence of the required commutative diagram is automatic. The result is proved.
\end{dem}

\section{(Re)visiting equivalences of categories in deformation theory}\label{sec deformation}

Let $k$ be a perfect field of characteristic $p$ and $X$ be a $k$-scheme. Recall that $\frob$ denotes the (absolute) Frobenius of $X$. 
We  present a functorial description of  three  problems in deformation theory, using the point of view of Witt modules:

\begin{enumerate}
\item Lift $X$ to a scheme $X_2$, flat over $\W_2(k)$.\vspace{0.3cm}\label{Pb1}
\item Lift $X$ to a scheme $X_2$, flat over $\W_2(k)$, \textit{together with its Frobenius morphism}. \label{Pb2}
\item Lift a given vector bundle $V/X$ to a $\W_2$-bundle $V_2$. \label{Pb3}
\end{enumerate}

In \eqref{Pb2}, one requires the existence of an endomorphism $F_2$ of $X_2$, whose mod $p$ reduction is the  Frobenius of $X$. This  is a very strong extra requirement. Since $k$ is perfect, such an $F_2$ is automatically compatible with the Frobenius of $\W_2(k)$.

Problems \eqref{Pb1} and \eqref{Pb2} have been the subject of sustained investigation by many authors---see, for instance, the seminal papers \cite{DI} and \cite{MS}. Most related  results below are  known. Nevertheless, our approach, via the extension $\mathcal E \W_2(X)$, is new.

Problem \eqref{Pb3} is relatively new.
 
\subsection{Recollections on deformation theory}
Let us first recall some well-known concepts and facts from deformation theory. For details, see \cite{Il}.

\begin{defi}Let $X$ be a scheme over a base $B$, and let $\mathcal M$ be a coherent $\mathcal O_X$-module. A square-zero extension of $X$ by $\mathcal  M$ (in the category of $B$-schemes) is the data of a $B$-scheme $Y$, together with a closed embedding
\[ i: X \to Y,\]
defined by an Ideal $\mathcal I \subset \mathcal O_Y $  of square zero, equipped with an isomorphism of $\mathcal O_X$-modules
\[f:  \mathcal I  \stackrel \sim \to \mathcal M. \]
Square-zero extensions of $X$ by $\mathcal M$ (in the category of $B$-schemes) form a category (where morphisms are isomorphisms), which we denote by $\mathbf{Exal_B}(X,\mathcal M).$ 
\end{defi}

\begin{rem}
The previous definition implicitly uses the fact that the $\mathcal O_Y$-module $\mathcal I$ is actually an  $\mathcal O_X$-module.

In short, a square-zero extension can be thought of as an extension \[ 0 \to \mathcal M \to \mathcal O_Y \stackrel \pi \to \mathcal O_X \to 0,\] in which $\mathcal M$ is a square-zero ideal, and $\mathcal \pi$ is a  homomorphism of $\mathcal O_B$-algebras. 
\end{rem}

\begin{defi}
Let $X$ be a scheme, and let $\M$ and $\mathcal{N}$ be $\mathcal O_X$-modules. Denote by $\EExt_{\mathcal O_X}^1( \mathcal{N}, \M)$ the category whose objects are  extensions of $\mathcal O_X$-modules
\[0 \to \M \to E \to \mathcal{N} \to 0, \]
 morphisms being morphisms of exact sequences, which are identity on $\M$ on $\mathcal{N}$. 
\end{defi}

\begin{rem}
Morphisms in $\EExt_{\mathcal O_X}^1( \mathcal{N}, \M)$ are isomorphisms, and the group of automorphisms of any object is naturally isomorphic to $\Hom_{\mathcal O_X}(\mathcal{N},\M).$
\end{rem}

\begin{prop}[\cite{Il}]\label{equcatdiff}
Let $X$ be a smooth $k$-variety, $\M$ be an $\mathcal O_X$-module and $\Omega^1_{X/k}$ be the sheaf of K\"ahler differentials of $X$ over $k$.\\ There is an equivalence of categories
\[\mathbf{Exal}_k(X,\M) \stackrel \sim \to \EExt_{\mathcal O_X}^1( \Omega^1_{X/k}, \M) \]
that is compatible with Baer sum.  
\end{prop}

\subsection{An equivalence of categories for Problem \eqref{Pb1}.\\}
To begin with, assume that $X/k$ is smooth. One can then replace ``flat'' by ``smooth'' in \eqref{Pb1}. By deformation theory,  there exists a natural class \[ \mathrm{Obs}(X_2) \in \Ext^2_{\mathcal O_X-\mathrm{Mod}}(\Omega^1_{X/k} , \mathcal O_X), \]  which vanishes if, and only if, Problem \eqref{Pb1} has a positive answer.

In what follows, we give a simple functorial interpretation of this class. Note that, in Proposition 2.2 of the recent article \cite{Yo}, a similar goal is achieved. The approach chosen there is quite different from ours: a main input in its proof is Proposition 1 of the Appendix of \cite{MS}, which is concerned with our Problem \eqref{Pb2}.

From now on, we remove the smoothness assumption on $X/k$.
\begin{lem}\label{flatlift}
Let $X$ be a $k$-scheme. The data of a lifting of $X/k$ to a scheme $X_2$, flat over $\W_2(k)$, is equivalent to that of a square-zero extension (of schemes over $\W_2(k)$)
\[(\mathcal E: 0 \to \mathcal O_X \to \mathcal O_Y \to \mathcal O_X \to 0) \in \mathbf{Exal}_{\W_2(k)}(X,\mathcal O_X),\]
such that $\kappa_{\mathcal E,p}=\Id$ (cf. Section \ref{sec kappa} for the definition of $\kappa$).
\end{lem}

\begin{dem}
The statement follows immediately from the following observation. If $M$ is a $\W_2(k)$-module, then it is flat if and only if it is free. Equivalently, for the exact sequence \[ \mathcal M: 0 \to pM \to M \to M/pM \to 0,\] the connecting arrow $\kappa_{\mathcal M,p}: M/pM \to pM$ is an isomorphism.
\end{dem}

\begin{rem}
In Lemma \ref{flatlift}, it is crucial to work in the category of $\W_2(k)$-schemes, and hence consider square-zero extensions in $\mathbf{Exal}_{\W_2(k)}(X,\mathcal O_X)$.
\end{rem}

\begin{defi}
Using the description of the previous Lemma, the liftings of $X/k$ to a scheme flat over $\W_2(k)$ form a category (with morphisms being isomorphisms). It is the full subcategory of $\mathbf{Exal}_{\W_2(k)}(X,\mathcal O_X)$ consisting of extensions having $\kappa= \Id$. Denote it by $\mathcal L_{\W_2(k)}(X),$ or simply by $\mathcal L_{2}(X)$.
\end{defi}
 
The notation $\mathcal{L}_2(X)$ is unambiguous, thanks to the following Lemma.
 
\begin{lem}
Let $X$ be a $k$-scheme. Then, the forgetful functor
\[ \mathcal L_{\W_2(k)}(X) \to \mathcal L_{\Z /p^2 \Z}(X)\]
is an equivalence of categories.
\end{lem}
 
\begin{dem}
It is enough to deal with the affine case $X=\Spec(A)$. Let $A_2$ be a lift of $A$ to a  flat (i.e. free) $(\Z/ p^2 \Z)$-algebra. There is a natural isomorphism of $k$-vector spaces
\begin{align*}
A_2/pA_2= A &\stackrel \sim \to pA_2,\\
\overline x &\mapsto px.
\end{align*}
Since $k$ is perfect, $A_2$ can be turned into a $\W_2(k)$-algebra in a unique way, by the classical formula
\begin{align*}
\W_2(k) &\to A_2\\
\tau(x^p) &\mapsto \tilde x^p,
\end{align*}
where $\tilde x\in A_2$ is any lift of $x \in k \subset A$. Furthermore, let $(E_i)_{i \in I} \in A_2^I$  be a family, such that $(e_i:=\overline{E_i})_{i \in I}$ is a $k$-basis of $A$. Let us check that $(E_i)_{i \in I}$ is a basis of  the $\W_2(k)$-module $A_2$. That $(E_i)_{i \in I}$ is generating, is a straightforward two-step d\'evissage. This is actually (and more generally) a  variant of Nakayama's Lemma, for possibly non-finitely generated modules over  Artinian local rings (here $\W_2(k)$). Suppose now there is a  relation
$$\sum_{i\in I} \tilde x_i E_i=0,$$
where $\tilde x_i \in A_2$ is non-zero for finitely many $i\in I$. Since $(e_i)$ is a $k$-basis of $A$ and $k$ is perfect, one may write  $\tilde x_i=p\tau(x_i)$ for $x_i \in k$. From the relation
$$p\sum \tau(x_i)E_i=0 \in A_2,$$
and via the isomorphism $p A_2 \stackrel \sim \to A,$ one gets  $$\sum x_i e_i=0 \in A,$$ hence $x_i=0$, implying that $\tilde x_i=0$. Thus, $A_2$ is a free $\W_2(k)$-module.\\ By uniqueness of the $\W_2(k)$-algebra structure, a lift of a morphism of $k$-algebras $A \to A'$ to a morphism of flat $(\Z/ p^2 \Z)$-algebras $A_2 \to A'_2$, is automatically $\W_2(k)$-linear. The Lemma is proved.
\end{dem}

Given this last result, the reader may thus do one of the following.

\begin{enumerate}
    \item Bluntly assume that $k= \F_p$ everywhere;\vspace{0.3cm}
    \item Go on with an arbitrary perfect field $k$, but not worry about checking $\W_2(k)$-linearity of homomorphisms.
\end{enumerate}

From now on, assume that $X$ is a \emph{reduced} $k$-scheme. The next definition is the key prerequisite for our equivalence of categories.

Recall that, for the natural exact sequence
\[\mathcal E \W_2(X): 0 \to \frob_*(\mathcal O_X) \to \W_2(\mathcal O_X) \to \mathcal O_X \to 0,\]
one has
\[\kappa_{\mathcal E \W_2(X),p}=\frob:\mathcal O_X \to \frob_*(\mathcal O_X),\] see the end of Section \ref{sec kappa}.

\begin{defi}\label{defext}
Denote by
\[ \mathcal E F(X): 0 \to \mathcal O_X \xrightarrow{\frob} \frob_*(\mathcal O_X) \stackrel d \to \frob_*(B^1_X) \to 0\]
the natural sequence of $\mathcal O_X$-modules, in which $\frob_*(B^1_X)$ is  the cokernel of Frobenius. This notation is in accordance with the usual one.

The pushforward of $\mathcal E \W_2(X)$ by $d$ is  a square-zero extension, denoted by
\[ \mathrm {C} \W_2(X): 0 \to \frob_*(B^1_X) \to \mathcal O_Y \to \mathcal O_X \to 0.\]
Clearly, it has $\kappa=0$. Hence, we have $p\mathcal O_Y=0$. In other words,
\[\mathrm {C} \W_2(X) \in \mathbf{Exal}_{k}(X,\frob_*(B^1_X)),\]
is a square-zero extension of $X$ in the category of $k$-schemes.

If $X/k$ is smooth, by Proposition \ref{equcatdiff} we get an extension of $\mathcal O_X$-modules
\[C \Omega(X): 0 \to \frob_*(B^1_X) \to E \to \Omega^1_{X/k} \to 0, \]
corresponding to $\mathrm {C} \W_2(X)$.
\end{defi}

Note that the extensions $\mathcal E F(X)$, $\mathrm {C} \W_2(X)$ and $C \Omega(X)$ \emph{naturally} depend on $X$.

\begin{rem}
The extension $C \Omega(X)$ is denoted this way as it is given by the Cartier operator, when $X/k$ is smooth. Note that $\mathrm {C} \W_2(X)$  (for Cartier-Witt) actually makes sense for any \textit{reduced} $X/k$.

\end{rem}

We  now state and prove the promised equivalence of categories. Recall that $X$ is any reduced $k$-scheme.

\begin{defi}
Let $d: \frob_*(\mathcal O_X) \to \frob_*(B^1_X)$ be as in Definition \ref{defext}. Denote by $\tilde  {\mathcal L} _2(X)$ the category whose objects are  pairs $(\mathcal E,f)$, consisting of a square-zero extension
\[(\mathcal E: 0 \to \frob_*(\mathcal O_X) \to \mathcal O_Y \to \mathcal O_X\to 0)  \in \mathbf{Exal}_{k}(X,\frob_*(\mathcal O_X)),\]
together with an isomorphism (in $\mathbf{Exal}_{k}(X,\frob_*(B^1_X))$)
\[f: d_*(\mathcal E) \stackrel \sim \to \mathrm {C} \W_2(X).\]
Morphisms in $\tilde {\mathcal L}_2(X)$ are (iso)morphisms commuting to the given isomorphisms.
\end{defi}

\begin{prop}\label{EquLiftW2} Let $X$ be a reduced scheme over $k$. There is an equivalence of categories
\[ \Phi: \mathcal L_2(X) \stackrel \sim \to  \tilde {\mathcal L}_2(X).\]
\end{prop}

\begin{dem}
Let
\[ (\mathcal E:  0 \to \mathcal O_X  \to \mathcal O_{X_2}\to \mathcal O_X \to 0) \in \mathcal L_2(X)\]
be a lift of $X$ to a scheme flat over $\W_2(k)$. One has $\kappa_{\mathcal E,p}= \Id$ (cf. Section \ref{sec kappa} for the definition of $\kappa$). The pushforward
\[  \frob_*(\mathcal E):  0 \to \frob_*(\mathcal O_X)  \to \tilde{\mathcal O}_{X_2} \to \mathcal O_X \to 0,\]
and
\[\mathcal E \W_2(X): 0 \to \frob_*(\mathcal O_X) \to \W_2(\mathcal O_X) \to \mathcal O_X \to 0,\]
are both square-zero extensions of $X$ by $\frob_*(\mathcal O_X)$, in the category of schemes over $\W_2(k)$. They both have $\kappa=\frob$, so that their Baer difference $ \mathcal E \W_2(X)-\frob_*(\mathcal E)$ is a square-zero extension
\[\mathcal F: 0 \to \frob_*(\mathcal O_X) \to \mathcal O_Y \to \mathcal O_X \to 0,\]
in the category of $k$-schemes (indeed, the $\kappa$ arrow commutes with Baer sum, so that $\mathcal F$ has $\kappa=0$). Because $d \circ \frob=0$, the extension $d_*(\frob_*(\mathcal E))$ has a natural trivialization, so that the square-zero extension $d_*(\mathcal F) $ is naturally isomorphic to $d_*(\mathcal E \W_2(X)) = \mathrm {C}\W_2(X)$. Denoting the natural isomorphism by $f$,  the assignment $$\mathcal E \to (\mathcal F,f)$$ gives rise to a functor $\Phi:\mathcal L_2(X)\to \tilde {\mathcal L}_2(X)$.

A quasi-inverse of $\Phi$ is obtained as follows. Pick an object $(\mathcal F,f)$ in $\tilde {\mathcal L}_2(X)$. View
\[\mathcal F: 0 \to \frob_*(\mathcal O_X) \to \mathcal O_Y \to \mathcal O_X\to 0,\]
as a square-zero extension in $\mathbf{Exal}_{\W_2(k)}(X,\frob_*(\mathcal O_X))$, having $\kappa=0$. Then, the Baer difference
\[\tilde {\mathcal E}:=(\mathcal E \W_2(X)- \mathcal F) \in \mathbf{Exal}_{\W_2(k)}(X,\frob_*(\mathcal O_X)),\]
has $\kappa= \frob$, and  $d_*(\tilde {\mathcal E})$ is equipped with the trivialization induced by $f$. Therefore, it naturally yields a square-zero extension
\[\mathcal E: 0 \to \mathcal O_X \to \mathcal O_{X_2} \to \mathcal O_X\to 0,\]
in $\mathbf{Exal}_{\W_2(k)}(X,\mathcal O_X)$, together with a natural isomorphism
\[\frob_*(\mathcal E) \stackrel \sim \to \tilde {\mathcal E}.\]
The extension $\mathcal E$ then has $\kappa=\Id$, hence belongs to  $\mathcal L_2(X)$. The assignment $(\mathcal F,f) \mapsto \mathcal E$ defines a functor $\Psi$ in the reverse direction which is the inverse of $\Phi$.
\end{dem}

We now present two consequences of the preceding result. The first one is the existence of liftings of Frobenius-split $k$-schemes, a question which has been studied by several authors. See for instance the recent preprint \cite[Theorem 4.4]{Yo}, where $X/k$ is assumed to be smooth, but where Frobenius-splitting is replaced by a weaker notion.

\begin{coro}[Lifting of Frobenius-split $k$-schemes]\label{LiftFrobSplit}
Let $X$ be a reduced $k$-scheme. Assume that $X$ is Frobenius-split, i.e. that the exact sequence of coherent $\mathcal O_X$-modules
\[ \mathcal E F(X): 0 \to \mathcal O_X \to \frob_*(\mathcal O_X) \to \frob_*(B^1_X) \to 0,\]
splits. Then, $X$ admits a lift to a scheme $X_2$, flat over $\W_2(k)$. More precisely, every splitting of $\mathcal E F(X)$ naturally determines such an $X_2$.
\end{coro}

\begin{dem}
Denote by $s:\frob_*(B^1_X)\to \frob_*(\mathcal O_X)$ the splitting of $\mathcal E F(X)$ and put $\mathcal F:=s_*(\mathrm {C} \W_2(X))$. It is then clear that $\mathcal F$ is a square-zero extension and that $d_*(\mathcal F)$ is isomorphic to $\mathrm {C} \W_2(X)$. Conclude by applying Proposition \ref{EquLiftW2}.
\end{dem}

For another proof of the next corollary, see \cite[Thm 9.5, Prop 1.1]{Sch}.

\begin{coro}\label{explicitI}
Assume that $X/k$ is smooth. Then giving a lift of $X$ to a scheme $X_2$, smooth over $\W_2(k)$, is equivalent to giving an extension of $\mathcal O_X$-modules
\[\mathcal E: 0 \to \frob_*(\mathcal O_X) \to E \to \Omega^1_{X/k} \to 0, \]
together with an isomorphism (in $\EExt^1_{\mathcal O_X}(\Omega^1_{X/k},\frob_*(B^1_X))$)
\[ d_*(\mathcal E) \stackrel \sim \to  C \Omega(X),\]
where $C \Omega(X) \in \EExt^1_{\mathcal O_X} (\Omega^1_{X/k}, \frob_*(B^1_X)) $ is the extension of Definition \ref{defext}.\\

In particular, such a lift $X_2$ exists if and only if the cup-product
\[\mathcal E F(X) \cup C \Omega(X) \in \EExt^2_{\mathcal O_X} (\Omega^1_{X/k}, \mathcal O_X),\]
vanishes.
\end{coro}

\begin{dem}
The first part of the Corollary is a translation of the equivalence of categories in  Proposition \ref{EquLiftW2}, using that of Proposition \ref{equcatdiff}. The second part  follows from the first, applying standard considerations in homological algebra.
\end{dem}

Note that the cup-product $\mathcal E F(X) \cup C \Omega(X) \in \EExt^2_{\mathcal O_X} (\Omega^1_{X/k}, \mathcal O_X)$ of the previous Corollary is the obstruction $\mathrm{Obs}(X_2)$, given by classical deformation theory.

\subsection{An equivalence of categories for Problem \eqref{Pb2}}
We  move on towards a functorial description of Problem \eqref{Pb2}: lifting $k$-schemes together with their Frobenius morphism.

\begin{prop}\label{LemLiftFrob}
Let $X$ be a reduced $k$-scheme. Let $X_2 \in \mathcal L_2(X)$ be a lift of $X$, viewed as a square-zero extension \[\mathcal E:  0 \to \mathcal O_X \to \mathcal O_{X_2} \to \mathcal O_X \to 0  \] in $\mathbf{Exal}_{\W_2(k)}(X,\mathcal O_X)$, having $\kappa=\Id$. Then the following hold.\\
\begin{enumerate}
\item Consider the pullback $\frob^*(\mathcal E) \in \mathbf{Exal}_{\W_2(k)}(X,\frob_*(\mathcal O_X))$. It is the extension whose middle term is the sheaf of $\W_2(k)$-algebras $\mathcal O_Y$, defined as the fibered product
\[\xymatrix{\mathcal O_Y \ar[r] \ar[d] & \mathcal O_{X} \ar[d]^{\frob} \\ \mathcal O_{X_2} \ar[r]  & \mathcal O_X.} \]
Then, this extension is naturally isomorphic to $\mathcal E \W_2(X)$.
\item The data of a lift of $X \xrightarrow{\frob} X$, to a morphism $X_2\xrightarrow{F_2} X_2$ is equivalent to the data of an isomorphism of square-zero extensions in $\mathbf{Exal}_{\W_2(k)}(X,\frob_*(\mathcal O_X))$,
\[ \frob_*(\mathcal E) \stackrel \sim \to \mathcal E \W_2(X).\]
\item Assume that $X$ is a smooth $k$-scheme. Denote by $T_X$ its tangent bundle. Then, there exists a natural class $$ \mathrm{obs}(F_2) \in H^1(X, \frob^*(T_X)),$$ which vanishes if and only if there exists an $F_2$, as in item (2).
\end{enumerate}
\end{prop}

\begin{dem}
By glueing, it is enough to deal with the case where $X= \Spec(A)$ is affine. Then $X_2= \Spec (A_2)$, for some $\W_2(k)$-algebra $A_2$, free as a  $\W_2(k)$-module.\\
There is a commutative diagram
\[ \xymatrix{ \mathcal E \W_2(A): 0 \ar[r]& \frob_*(A) \ar[r] \ar[d]^{\Id} & \W_2(A) \ar[r] \ar[d]^{f_2} & A \ar[r] \ar[d]^{\frob} & 0 \\\mbox{\hspace{1.1cm}}\mathcal E:  0 \ar[r]& A \ar[r] & A_2  \ar[r]  & A \ar[r]  & 0,} \]
where $f_2$ is the natural ring homomorphism. The slightly abusive notation $$\Id: \frob_*(A) \to A$$ makes sense, remembering that  $\frob_*(A)=A$, as Abelian groups. The existence of this diagram proves item (1). For item (2), observe that a lift of $A \xrightarrow{\frob} A$, to a ring homomorphism $A_2 \xrightarrow{F_2} A_2$, is equivalent to a commutative diagram
\[ \xymatrix{ \mathcal E : 0 \ar[r]& A \ar[r] \ar[d]^{\frob} &  A_2 \ar[r] \ar[d]^{F_2} & A \ar[r] \ar[d]^{\frob} & 0 \\\mathcal E:  0 \ar[r]& A \ar[r] & A_2  \ar[r]  & A \ar[r]  & 0,} \]
in other words, an isomorphism  (in $\mathbf{Exal}_{\W_2(k)}(A,\frob_*(A))$)
\[ \frob_*(\mathcal E) \stackrel \sim \to \frob^*(\mathcal E).\]
 By item (1), the right side is naturally isomorphic to $\mathcal E \W_2(A)$, whence the claim. To prove item (3), form the Baer difference $$\Delta:= (\frob_*(\mathcal E)-\mathcal E \W_2(A)): 0 \to \frob_*(A) \to D \to A \to 0.   $$ Since formation of the connecting arrow $\kappa$ commutes to Baer sum,  one  has $\kappa_\Delta=\frob_A-\frob_A=0$. Thus, $D$ is a $k$-algebra, i.e. $\Delta$ belongs to $\mathbf{Exal}_k(A,\frob_*(A))$. Since $A/k$ is smooth, it is given by an extension of $A$-modules $$ 0 \to \frob_*(A) \to * \to \Omega^1(A/k) \to 0, $$ see Proposition \ref{equcatdiff}. Recalling that the vector bundles $T_X$ and $\Omega^1(X/k)$ are dual to each other, one gets $$\Ext^1_{A-\mathrm{Mod}}(\Omega^1(A/k),\frob_*(A))=H^1(X, \Omega^1(X/k)^\vee  \otimes_{\mathcal O_X} \frob_*(\mathcal O_X))$$ $$=H^1(X,  T_X \otimes_{\mathcal O_X} \frob_*(\mathcal O_X))=H^1(X,  \frob^*(T_X) ),$$ where the last equality uses  the projection formula. The claim follows.
\end{dem}

\begin{prop}[Obstruction for Problem \eqref{Pb2}]\label{EquLiftFrobW2} Let $X/k$ be a reduced scheme. The data of a lift of $X$ to a scheme $X_2$,  flat over $\W_2(k),$ together with a lift $\frob_2 :X_2 \to X_2$ of the Frobenius of $X$, is equivalent to that of a splitting of the square-zero extension $$C \W_2(X) \in \EExal^1_{k}(X, \frob_*(B^1_X) ).$$ If $X/k$ is smooth, this is equivalent to the data of a splitting of  the extension  $$C \Omega(X) \in \EExt^1_{\mathcal O_X}(\Omega^1_{X/k}, \frob_*(B^1_X) ).$$
\end{prop}

\begin{dem}
Use the equivalence of categories provided in  Proposition \ref{EquLiftW2}.  Keeping the notation of its proof, we see by Proposition \ref{LemLiftFrob}  that the data of a  lift of $X$, flat over $\W_2(k)$, together with its Frobenius, amounts to specifying an isomorphism $ \frob_*(\mathcal E) \stackrel \sim \to \mathcal E \W_2(X);$ that is, a splitting of $C \W_2(X)$. To prove the last assertion (when $X/k$ is smooth), remember that $C \Omega(X)$ then corresponds to $C \W_2(X)$, through the equivalence of Proposition \ref{equcatdiff}.
\end{dem}

\subsection{An equivalence of categories for Problem \eqref{Pb3}}
The approach taken here is, \textit{mutatis mutandis}, the same as that used to tackle Problem \eqref{Pb1}. Some proofs are thus left to the reader. Symmetric and divided powers of modules are freely used below. These are polynomial functors, characterized by a universal property (see \cite{Fe} for details). \\

Let $X$ be a scheme over $k$. Let $V$ be a vector bundle over $X$. Denote by $$f: P:=\P_X(V) \to X$$ the projective bundle of $V$. Denote by \[\mathrm{ad}: V \to \frob_*(\frob^*(V)) \] the $\mathcal O_X$-linear (first) adjunction morphism. Recall that, if  $\M$ a quasi-coherent $\mathcal O_X$-module, we have an adjunction isomorphism \[\EExt^1_{\mathcal O_X}(\frob^*(V),\M) \stackrel \sim \to  \EExt^1_{\mathcal O_X}(V ,\frob_*(\M)),\]  given by applying the exact functor $\frob_*$, followed by pulling back by $\mathrm{ad}$. To get its inverse, apply $\frob^*$, and push forward by the (second) $\mathcal O_X$-linear adjunction $$\frob^*(\frob_*(V)) \to V. $$ 

\begin{defi}
There are morphisms (of vector bundles over $X$)
\begin{align*}
    \Ver_V: \frob^*(V) &\to \Sym_{\mathcal O_S}^p(V)\\
    v \otimes 1 &\longmapsto v^p
\end{align*}
and
\begin{align*}
    \frob_V: \Gamma_{\mathcal O_S}^p(V) &\to\frob^*(V)\\
    [v]_p   &\longmapsto v \otimes 1
\end{align*}
called the \emph{Verschiebung} and the \emph{Frobenius} of $V$, which fit into exact sequences  \[ \mathcal E \Ver(V): 0 \to \frob^*(V) \xrightarrow{\Ver_V} \Sym_{\mathcal O_X}^p(V) \stackrel {q_V}\to   \overline \Sym_{\mathcal O_X}^p(V) \to 0\] and \[ \mathcal E \frob(V):0 \to  \overline  \Gamma_{\mathcal O_X}^p(V)\to  \Gamma_{\mathcal O_X}^p(V) \xrightarrow{\frob_V} \frob^*(V) \to 0,\]  where $ \overline \Sym_{\mathcal O_X}^p(V)$  (resp. $\overline  \Gamma_{\mathcal O_X}^p(V)$) is defined as as the cokernel of $\Ver_V$ (resp. kernel of $\frob_V$).\\
\end{defi}

\begin{rem}
The dual of the exact sequence (of vector bundles over $X$) $\mathcal E \Ver(V^\vee)$ is naturally isomorphic to $\mathcal E \frob(V).$
\end{rem}

\begin{lem}
There is a natural isomorphism (of vector bundles over $X$) \[\Phi_V: \overline \Sym_{\mathcal O_X}^p(V)  \stackrel \sim \to \overline  \Gamma_{\mathcal O_X}^p(V).\]
\end{lem}

In what follows, we may tacitly use this result to identify these vector bundles.

\begin{dem}
Consider the natural homomorphism
\[ \alpha_V: \Sym_{\mathcal O_X}^p(V)  \to   \Gamma_{\mathcal O_X}^p(V),\]
defined on sections by the formula
\[ v_1 v_2 \ldots v_p \mapsto [v_1]_1 \ldots [v_p]_1.\]
It takes values in $\overline  \Gamma_{\mathcal O_X}^p(V)$, and vanishes on $\Im(\Ver_V)$ because of the identity $$[v]_1^p=p! [v]_p=0, $$  which follows  by induction from the rule $$[x]_i [x]_j= {{i+j} \choose i}[x]_{i+j},$$ which is part of the definition of divided powers (see \cite{Fe}).\\

The resulting homomorphism \[ \overline \Sym_{\mathcal O_X}^p(V)  \to \overline  \Gamma_{\mathcal O_X}^p(V)\] is an isomorphism. To check this, one can assume that $V=\mathcal O_X^d$ is the trivial rank $d$ vector bundle, and that $X=\Spec(A)$ is affine. The rest of the verification is then straightforward. Indeed,  the description of $\Sym^p$ and $\Gamma^p$ as polynomial functors then presents both $ \overline   \Sym_{\mathcal O_X}^p(V)$ and $ \overline   \Gamma_{\mathcal O_X}^p(V)$ as trivial vector bundles, with respective canonical basis $$e_1 ^{a_1} \ldots e_p^{a_p}$$ and  $$[e_1]_{a_1} \ldots [e_p]_{a_p},$$ both indexed by proper partitions \[ a_1+ \ldots +a_p =p.\] Here ``proper'' means that at least two $a_i$'are non-zero. Since $(p-1)!$ is invertible, $\alpha_V$ maps each $e_1 ^{a_1} \ldots e_p^{a_p}$ to an invertible multiple of $[e_1]_{a_1} \ldots [e_p]_{a_p}$ and the Lemma is proved.
\end{dem}

\begin{rem}
By adjunction, we have a natural isomorphism \[\EExt^1_{\mathcal O_X}(\frob^*(V) ,\overline  \Gamma_{\mathcal O_X}^p(V)) \stackrel \sim \to \EExt^1_{\mathcal O_X}(V ,\frob_*(\overline  \Gamma_{\mathcal O_X}^p(V))), \] through which $\mathcal E \frob(V)$ corresponds to an extension  \[ \overline {\mathcal E } \frob(V):0 \to  \frob_*(  \overline  \Gamma_{\mathcal O_X}^p(V))\to  \Phi(V)  \to  V \to 0.\] 
\end{rem}

To get a mod $p^2$ avatar of  $\overline {\mathcal E } \frob(V)$, one uses  Teichm\"uller lifts of line bundles.\\
Recall (Section \ref{sec TeichLift}) the natural exact sequence of Witt modules on $P:=\P(V)$ \[0 \to \frob_*(\mathcal O_P(p)) \to \W_2(\mathcal O_P(1)) \to \mathcal O_P(1) \to 0. \]  Applying $f_*$, one gets an exact sequence of Witt modules on $X$   \[0 \to \frob_*(\Sym_{\mathcal O_X}^p(V)) \to f_*( \W_2(\mathcal O_P(1)) ) \stackrel {\Phi_V }\to  V \to 0. \] Surjectivity of the last arrow follows from the (effect on global sections of the)  Teichm\"uller section $\tau_{\mathcal O_P(1)}$; see Section \ref{TeichWitt}.

\begin{defi}
The exact sequence of $\W_2(\mathcal O_X$)-modules on $X$
\[0 \to \frob_*(\Sym_{\mathcal O_X}^p(V))\stackrel {i_V} \to f_*( \W_2(\mathcal O_P(1)) ) \stackrel {\Phi_V }\to  V \to 0,\]
is denoted by $\overline {\mathcal E} \W_2(V)$. Denote by $\tau_V$ the natural (sheaf-theoretic) section  of $\Phi_V$ induced by the  Teichm\"uller section $\tau_{\mathcal O_P(1)}$.
\end{defi}

\begin{rem}
Assume that $X=\Spec(A)$ is affine. Denote by $B$ the symmetric algebra  $\Sym_A(V)=\bigoplus_{i=0}^\infty \Sym^i_A(V)$. One can also obtain $\overline {\mathcal E} \W_2(V)$ from the exact sequence \[0 \to \frob_*(B) \to \W_2(B) \to B \to 0, \] by pulling it back by the inclusion $V \to B$, and pushing it forward by the projection $  \frob_*(B) \to  \frob_*(\Sym_A^p(V))$.
\end{rem}

\begin{lem}
One has the formula
\[\tau_V(x+y)=\tau_V(x)+\tau_V(y)+i_V\left(\sum_{i=1}^{p-1} \frac 1 p \binom{p}{i}  x^i y^{p-i}\right).\]
\end{lem}

\begin{dem}
Follows from th fact that the same formula defines addition on $\W_2$.
\end{dem}

\begin{lem}
The extension $\overline {\mathcal E} \W_2(V)$ has $\kappa:=\kappa_{\overline{\mathcal E} \W_2(V),p}$ given by the map adjoint to $\Ver_V$. Concretely, it is given by
\[\fonctionnoname{V}{\frob_*(\Sym_{\mathcal O_X}^p(V)).}{x}{x^p}\]
There is a natural commutative diagram of $\W_2(\mathcal O_S)$-modules
\[\xymatrix@C=2em{ 0 \ar[r] & \frob_*(\Sym_{\mathcal O_X}^p(V)) \ar[r] \ar[d]^{\frac 1 {(p-1)!}\frob_*(q_V)}  & f_*( \W_2(\mathcal O_P(1)) ) \ar[r] \ar[d]^F & V \ar[r] \ar[d]^{\mathrm{ad}} & 0 \\ 0 \ar[r] & \frob_*( \overline \Sym_{\mathcal O_X}^p(V)) \ar[r] & \frob_*(\Gamma^p_{\mathcal O_X}(V))  \ar[r]  &  \frob_*(\frob^*(V)) \ar[r] & 0 ,} \]
where the upper line is $\overline {\mathcal E }\W_2(V)$, the lower line is $\frob_*({\mathcal E } \frob(V)),$ and  where $F$ is defined (on sections) by the formula  $$F(\tau_V(x))=  [x]_p,$$ for $x \in V.$

Consequently, $(\frob_*(q_V))_*(\overline {\mathcal E }\W_2(V))$  is naturally isomorphic to $-\overline {\mathcal E } \frob(V)$.
\end{lem}

\begin{dem}
Can assume  $X=\Spec(A)$ is affine. The key point here is  that the formula giving $F$ makes sense, and indeed defines a homomorphism of Witt modules. This follows from the previous Lemma, combined to the formula, in $\Gamma^p_A(V)$, \[ [x+y]_p=[x]_p+[y]_p+\sum_{i=1}^{p-1}  [ x]_i [y]_{p-i}.\] \end{dem} 

\begin{thm}[An equivalence of categories for Problem \eqref{Pb3}]\label{EquW2}
To give a lift of $V$ to a $\W_2$-bundle $V_2$ on $X$ is equivalent to give an extension
\[ \mathcal F \in \EExt^1_{\mathcal O_X}(\frob^*(V) ,  \Sym_{\mathcal O_X}^p(V)),\]
together with an isomorphism
\[ (q_V)_*(\mathcal F) \stackrel \sim \to {\mathcal E } \frob(V),\]
in $\EExt^1_{\mathcal O_X}(\frob^*(V) ,  \overline \Sym_{\mathcal O_X}^p(V))$.
\end{thm}

\begin{dem}
Similar to that of Proposition \ref{EquLiftW2}.
\end{dem}

In purely cohomological terms, one immediately deduces the following corollary.

\begin{coro}\label{explicitII}
The obstruction $\mathrm{Obs}(V_2)$ to lifting $V$ to a $\W_2$-bundle on $X$ is the element of $\Ext^2_{\mathcal O_X}(\frob^*(V),\frob^*(V) )$ represented by the $2$-extension \[ 0 \to \frob^*(V)  \to \Sym^p_{\mathcal O_X} (V) \to \Gamma^p_{\mathcal O_X} (V) \to   \frob^*(V) \to 0,\] defined as the cup-product of $\mathcal E \frob(V)$ and $\mathcal E \Ver(V)$. 
\end{coro}

\subsection{Relating Problems \eqref{Pb2} and \eqref{Pb3}.}

\begin{prop}\label{exoconnection} 
Let $X$ be a smooth variety, over a perfect field $k$, $\carac(k)=p$. Denote by $T_X$ the tangent bundle of $X/k$. Let $V$ be a vector bundle on $X$.\\ Consider the following assertions.
\begin{enumerate}[(a)]
    \item Lift the variety $X$ to a smooth scheme over $\W_2(k)$, together with its Frobenius.\\ \label{(a)}
    \item Lift the variety $X$  to a smooth scheme  $X_2/\W_2(k)$, in such a way that $\frob^*(V)$ lifts, to a vector bundle on $X_2$.\\ \label{(b)}
    \item Lift the vector bundle $V/X$ to a $\W_2$-bundle on $X$. \label{(c)}
\end{enumerate}
Then \eqref{(a)} implies \eqref{(b)} for $V=T_X$, and \eqref{(b)} implies \eqref{(c)}.
\end{prop}

\begin{dem}
The first implication is straightforward. Indeed, let $X_2 \to \Spec(\W_2(k))$ be a smooth morphism, lifting $X \to \Spec(k)$. Denote by $V_{X_2}$  its tangent bundle. It is a vector bundle over $X_2$. Let $ X_2 \xrightarrow{F_2} X_2$ be a lift of $ X \xrightarrow{\frob} X$. Then, $F_2^*(V_{X_2})$ is a lift of $\frob^*(V)$, to a vector bundle on $X_2$.\\
Let us prove  \eqref{(b)}$\Longrightarrow$\eqref{(c)}. Let $V_{X_2}^{[1]}$ be a lift of $\frob^*(V)$, to a vector bundle over $X_2$. There is the natural reduction sequence of $\mathcal O_{X_2}$-modules $$0 \to \frob^*(V) \to V_{X_2}^{[1]} \to \frob^*(V) \to 0.$$ Recall that, at the level of topological spaces, $X_2=X$ and   $ X \xrightarrow{\frob} X$ is the identity. Consider the extension above, merely as  an extension of Zariski sheaves on $X$. One may then replace $\frob^*(V)$ by $\frob_*(\frob^*(V))$. Form the pull-back diagram of Zariski sheaves on $X$\[ \xymatrix{ \mathcal E V_2: 0 \ar[r]& \frob_*(\frob^*(V)) \ar[r] \ar@{=}[d] & V_2 \ar[r] \ar[d] & V \ar[r] \ar[d]^{\frob_V} & 0 \\  0 \ar[r]&  \frob_*(\frob^*(V))\ar[r] & V_{X_2}^{[1]} \ar[r]  &  \frob_*(\frob^*(V)) \ar[r]  & 0,} \]  where $\frob_V$ is the natural $\mathcal O_X$-linear arrow. Using the first commutative diagram of the proof of Proposition \ref{LemLiftFrob}, one sees that $V_2$ has a natural structure of a $\W_2(\mathcal O_X)$-module,  for which $\mathcal E V_2$  is an extension of  $\W_2(\mathcal O_X)$-modules.  It is straightforward to check, that its associated connecting arrow $V \xrightarrow{\kappa} \frob_*(\frob^*(V))$ equals $\frob_V$. Hence, $V_2$ is a lift of $V$ to a $\W_2$-bundle, and $\mathcal E V_2$  is the reduction sequence of $V_2$.
\end{dem}

\begin{rem}
Assume there exists a lift $X_2$, as in \eqref{(b)} above. Then liftability of $V$ to a  $\W_2$-bundle on $X$, implies that of $\frob^*(V)$ to a vector bundle on $X_2$. This is a particular case of Lemma \ref{LemLiftUniv}.
\end{rem}

\section{Grassmannians whose tautological bundle does not lift}

In this section,   Grassmannian varieties $\Gr(m,n)$ are considered over a base of characteristic $p$. Our goal  is the next Theorem, proved by  applying the equivalences of categories offered in the previous section, to tautological bundles on $\Gr(m,n)$.

\begin{thm}\label{nonlift}
Let $m$ and $n$ be two integers, with $2 \leq m \leq n-2$. \vspace{0.2cm}\\
Then, the tautological vector bundle $V$ of $\Gr(m,n)$ does not lift to a $\W_2$-bundle.\\Neither do its Frobenius twists $V^{(s)}$, for all $s \geq 0$.
\end{thm}
\begin{rem}
In  the Proposition,  one can assume without loss of generality that the base is a perfect field $k$. Using that formation of coherent cohomology of varieties commutes with change of the base field $k$, one could actually reduce to $k=\F_p$.
\end{rem}

\subsection{First proof of Theorem \ref{nonlift}}\quad\vspace{0.3cm}\\
We give here a first proof of Theorem \ref{nonlift} that uses Proposition \ref{EquW2} and some explicit computations that we carry in Section \ref{sec computations}. We only proof the case where $s=0$, the general case being the same.

Let $k$ be a perfect field of characeristic $p$. Let $E$ be a $k$-vector space of dimension $n$ and let $X:=\Gr(m,E)$ be the corresponding Grasmannian, which parametrizes $m$-dimensional subspaces of $E$. Denote by $f:X\to \Spec(k)$ its
structure morphism and by $V$ its tautological bundle. There is an exact sequence
\[ 0 \to V \to f^*(E) \to W \to 0,\]
whose cokernel is a vector bundle $W$ on $X$. One knows (see Section \ref{sec TautoLift}) that $V$ admits a lifting tower, in the (dual) cases $m=1$ or $m=n-1$.\vspace{0,3cm}\\
Assume now that $2 \leq m \leq n-2$. The goal is to show that $V$ does not lift to a $\W_2$-bundle on $X$. Consider the vector bundle extension \[ \mathcal E \frob(V):0 \to  \overline  \Sym_{\mathcal O_X}^p(V)\to  \Gamma_{\mathcal O_X}^p(V) \xrightarrow{\frob_V} \frob^*(V) \to 0.\] Assume that $V$ can be lifted to a $\W_2$-bundle. By Proposition \ref{EquW2}, this means that $ \mathcal E \frob(V)$  admits a lift to an extension \[\mathcal F: 0 \to    \Sym_{\mathcal O_X}^p(V)\to  F  \to   \frob^*(V) \to 0.\]   By item 2) of  Lemma \ref{calcul4} (for $r=1$), $\mathcal F$ would be split, hence so would $\mathcal E \frob(V)$. But this contradicts item 4) of the same Lemma (again for $r=1$).

\subsection{A second proof for Theorem \ref{nonlift}}\quad\vspace{0.3cm}\\
We give now a second proof that promotes the use of Teichm\"uller lifts of line bundles. Though  more elementary than the first one, it relies on the specific fact that our base $X$ is a projective homogeneous space of a reductive algebraic group. Here again, we assume for simplicity that $s=0$, the proof of the general case being the same, and we use computations from Section \ref{sec computations}.

Keep notation of the preceding section.  In particular, one has $2 \leq m \leq n-2$ and $X=\Gr(m,E)$ with structure morphism $f: X \to \Spec(k)$.

Assume that $V$ admits a lift to a $\W_2$-bundle on $X$. Consider the  arrow
$$ h: \Fl(E)=\Fl(1,2,\ldots, n-1, E) \to X,$$
where $\Fl(E)$ denotes the complete flag scheme of $E$ (see Definition \ref{defiflag}). Denote the tautological filtration on $\Fl(E)$ by \[ 0 \subset \mathcal V_1 \subset  \mathcal V_2 \subset \ldots \subset   \mathcal V_n=E, \]and its graded pieces by $\mathcal L_i:= \mathcal V_i /\mathcal V_{i-1}$. Thus, over $\Fl(E)$, $h^*(V)=\mathcal V_m$, so that  $\mathcal V_m$ admits a lift, to a $\W_2$-bundle $\mathcal V_{m,2}$ over $\Fl(E)$.

Assume that $m\geq 3$. We claim that the natural quotient
\[ \pi_m: \mathcal V_m \to \mathcal L_m,\]
lifts to a surjection of $\W_2$-bundles
\[ \pi_{m,2}: \mathcal V_{m,2} \to \W_2(\mathcal L_m).\]
To prove it, one proceeds as in the proof of Theorem \ref{PropTautoLift}. The space of such lifts is a torsor under the vector bundle $\mathcal V_m^{\vee(1) }\otimes \mathcal L_m^{(1)}$. By item (2) of Lemma \ref{calcul2}, this torsor is trivial. The claim follows.\vspace{0.3cm}\\
The kernel of $\pi_{m,2}$ is then a lift of $\mathcal V_{m-1}$ to a $\W_2$-bundle $\mathcal V_{m-1,2}$. By descending induction on $m$, one infers that  $\mathcal V_{2}$ lifts, to a $\W_2$-bundle $\mathcal V_{2,2}$ over $\Fl(E)$, and that the arrow \[ \pi_2: \mathcal V_2 \to \mathcal L_2\] lifts to a surjection of $\W_2$-bundles \[ \pi_{2,2}: \mathcal V_{2,2} \to \W_2(\mathcal L_2),\] whose kernel $\mathcal L_{1,2}$ is a $\W_2$-line bundle, lifting  $\mathcal L_{1}$.  By Proposition \ref{LiftToW2}, the space of such lifts is a principal homogeneous space of $H^1(\Fl(E),\mathcal O_{\Fl(E)})$, pointed by $\W_2(\mathcal L_1)$. As $H^1(\Fl(E),\mathcal O_{\Fl(E)})=0,$  one sees that $\mathcal L_{1,2}$ is isomorphic to $\W_2(\mathcal L_1)$.\\ Altogether, we have built an extension of $\W_2$-bundles \[\mathcal E_{1,2,2}:= 0 \to \W_2(\mathcal L_1)\to \mathcal V_{2,2} \to \W_2(\mathcal L_2) \to 0,\] lifting the tautological extension \[\mathcal E_{1,2} : 0 \to \mathcal L_{1} \to \mathcal V_{2} \to \mathcal L_2 \to 0.\] The natural exact sequence  of Witt modules \[0 \to (\mathcal L_2^\vee \otimes \mathcal L_1) ^{(1)} \to \W_2(\mathcal L_2^\vee \otimes \mathcal L_1)  \to \mathcal L_2^\vee \otimes \mathcal L_1 \to 0\]  admits a sheaf-theoretic section:  the Teichm\"uller section $\tau$. Thus, the sequence \[0 \to H^1(X, (\mathcal L_2^\vee \otimes \mathcal L_1) ^{(1)}) \stackrel \iota \to H^1(X, \W_2(\mathcal L_2^\vee \otimes \mathcal L_1) ) \to H^1(X, (\mathcal L_2^\vee \otimes \mathcal L_1) ^{(1)}) \] is exact. Therefore, the set of isomorphism classes of lifts of $\mathcal E_{1,2} $  to an extension $\mathcal E_{1,2,2}$ as above, is a principal homogeneous space of $$\Ext^1(\mathcal L_2^{(1)}, \mathcal L_1^{(1)})=H^1(X, (\mathcal L_2^\vee \otimes \mathcal L_1) ^{(1)}).$$ By point (3) of Lemma \ref{calcul2}, this is a one-dimensional $k$-vector space, with generator (the class of) $\mathcal E_{1,2}^{(1)}$. Note that the computation  (3) of Lemma \ref{calcul2} remains valid over an arbitrary base ring of characteristic $p$. Hence, the result above remains valid after changing the base, from $k$ to an arbitrary $k$-algebra.

Set $G:=\GL_k( E)$ for the linear algebraic $k$-group of  linear automorphisms of $E$. Consider $\Fl(E)$ as a projective homogenous space of $G$, over $k$. For this structure, the tautological vector bundles $\mathcal V_i$, as well as the line bundles   $\mathcal L_i$, are naturally $G$-linearized. By functoriality of the Teichm\"uller lift of line bundles, the $\W_2$-bundles $\W_2(\mathcal L_i)$, over $\Fl(E)$, are $G$-linearized as well. Thus, $\mathcal E_{1,2,2}$ is an extension between the $G$-linearized $\W_2$-bundles $\W_2(\mathcal L_2)$ and $\W_2(\mathcal L_1)$. However, so far, its middle term $\mathcal V_{2,2} $ need not admit a $G$-linearization. We are going to show that this can in fact be done, in the strongest possible sense.

Denote by $k':=k[G]$ the function ring of $G$, and by $g \in G(k')$ the point corresponding to the identity of $G$. Change the base ring, from $k$ to $k'$. Set $$\Fl(E)':=\Fl(E) \times_k k', $$and similarly for other objects, and work over $k'$.   Consider the base-changed extension
\[\mathcal E_{1,2,2}':= 0 \to \W_2(\mathcal L_1')\to \mathcal V'_{2,2} \to \W_2(\mathcal L_2') \to 0,\]
of vector bundles over $\Fl(E)'$. Set $$\fonction{g^*}{\Fl(E)'}{\Fl(E)'}{\nabla}{(g^{-1}).\nabla}$$
where the presence of an exponent $-1$ guarantees that pulling back by (scheme-theoretic) points of $G$, provides a \textit{left  action} of $G$, on $\W_2$-bundles.

Consider the extension

{\centering\noindent\makebox[355pt]{$g^*(\mathcal E_{1,2,2}'):= 0 \to \W_2(\mathcal L_1') \stackrel \sim \to g^*(\W_2(\mathcal L_1'))\to g^*(\mathcal V'_{2,2}) \to g^*(\W_2(\mathcal L_2') ) \stackrel \sim \to \W_2(\mathcal L_2') \to 0,$}}

where isomorphisms are given by the natural $G$-linearizations. The set of isomorphism classes of lifts of $\mathcal E_{1,2}' $ to $\EExt^1_{G,2}(\W_2(\mathcal L_2'),\W_2(\mathcal L_1'))$  is a principal homogeneous space of $k'$.  Precisely, there exists a unique element $\lambda \in k'$, such that
\[g^*(\mathcal E_{1,2,2}')-\mathcal E_{1,2,2}' \simeq \lambda \iota(\mathcal E_{1,2}^{'(1)} ) \in \EExt^1_{G,2}(\W_2(\mathcal L_2'),\W_2(\mathcal L_1')). \]
Considered as a morphism of $k$-varieties,
\[ \lambda: G \to \G_a\]
is then a Hochschild $1$-cocycle, i.e. a group homomorphism. Since
$$\Hom(G,\G_a)=0,$$
it is trivial. Hence, the extensions (of $\W_2$-bundles over $\Fl(E)'$) $g^*(\mathcal E_{1,2,2}')$ and $\mathcal E_{1,2,2}'$ are isomorphic. Let
$$\phi_g: g^*(\mathcal E_{1,2,2}') \to \mathcal E_{1,2,2}'$$
be an isomorphism of extensions. Denote the coordinates of $G \times_k G$ by $(g_1,g_2)$, set $k'':=k[G \times_k G]=k' \otimes_k k'$ and $\Fl(E)''=\Fl(E) \times_k k''$. Working over $k''$, the formula \[\phi_{g_1} \circ (g_1)^* (\phi_{g_2}) \circ \phi_{g_1 g_2}^{-1} \] defines an  automorphism of the extension  $\mathcal E_{1,2,2}''$; that is, an element of  \[ \Hom_{\Fl(E)''}(\W_2(\mathcal L_2''), \W_2(\mathcal L_1'')),\] which is trivial by a straightforward two-step d\'evissage, using point (4) of Lemma \ref{calcul2} below. Thus, the cocycle condition \[\phi_{g_1} \circ (g_1)^* (\phi_{g_2}) =\phi_{g_1 g_2} \] identically holds. In other words,  $\phi_g$ provides a $G$-linearization of   $\mathcal E_{1,2,2}$. \\Choose a  $k$-basis of $E$, i.e. a direct sum decomposition of $k$-vector spaces \[ E=L_1  \bigoplus L_2  \bigoplus \ldots \bigoplus L_n\] with $L_i=k$ for all $i$. It determines an embbedding of $k$-groups \[ \G_a \subset \GL_k(E) \simeq \GL_n,\] identifying  $\G_a $ to  the subgroup consisting of strictly upper triangular matrices \[\begin{pmatrix} 1 & x& 0 & \ldots & 0 \\  0 & 1 & 0 & \ldots &  0 \\   0 & 0 & 1 &  \ldots &  0 \\ \vdots  &\vdots  & \vdots & \ddots &  \vdots  \\ 0 & 0 & 0 &  \ldots &  1\\
 \end{pmatrix}. \] The filtration associated to the chosen  basis  gives a $\G_a$-invariant $k$-rational point $$\nabla \in \Fl(E)(k).$$  To conclude the proof, let us show how specialising at $\nabla$ yields a contradiction.\\ The fiber of  $\mathcal E_{1,2}$, at $\nabla$, is the natural (non-split) extension of $\G_a$-modules over $k$, $$E_{1,2}: 0 \to L_1 \to L_1 \oplus L_2 \to   L_2 \to 0,$$ where $\G_a$ acts as the vector group of the one-dimensional $k$-vector space $\Hom(L_2,L_1)$. Similarly, the fiber of  $\mathcal E_{1,2,2}$ at $\nabla$ is a $\G_a$-linearized lift of $E_{1,2}$, to a $\G_a$-linearised extension of $\W_2(k)$-modules $$E_{1,2,2}: 0 \to \W_2(L_1) \to \W_2(L_1) \oplus \W_2(L_2) \to   \W_2(L_2) \to 0,$$ with trivial action of $\G_a$ on both graded pieces.  This action is simply given by a homomorphism of linear algebraic $k$-groups $$\alpha_2: \G_a \to \W_2(L_2 ^\vee \otimes L_1) \simeq R_{\W_2(k)/k}(\G_{a,\W_2(k)}).$$ Here $\W_2(L_2 ^\vee \otimes L_1)$  is considered as a two-dimensional group scheme over $k$, isomorphic to the Greenberg transfer (see \cite{BGA}) of  $\G_a$, from $\W_2(k)$ to $k$.\\ That  $E_{1,2,2}$ lifts $E_{1,2}$, means that the composite arrow of algebraic $k$-groups $$ \G_a \xrightarrow{\alpha_2} R_{\W_2(k)/k}(\G_{a,\W_2(k)}) \xrightarrow{\rho} \G_a$$ is the identity. (Here $\rho$ is the natural reduction homomorphism.)\\
 Taking $k$-points, one gets a factorisation of $\Id_k$, as a composite of homomorphism of additive groups $$ k \xrightarrow{\alpha_2(k)} \W_2(k) \to k,$$ where the right arrow is the reduction homomorphism. Such a splitting does not exist, because $p\neq 0 \in \W_2(k)$. The proof is complete.\\

\section{Application: non-liftability of some projective bundles}\label{Appli2}

Theorem \ref{nonlift} is close to Theorem 6.5 of \cite{Z}. The connection is made explicit below. \\
The result  is valid when $\F_p$ is replaced by an arbitrary field $k$ of characteristic $p$, with the same proof.
\begin{thm} \label{Zda}(see \cite{Z}, Theorem 6.5)\\
Let $m,n$ be integers, with $2 \leq m \leq n-2$.  Denote by $V$ the tautological vector bundle on the Grassmannian $$X:=\Gr(m,n) \to \Spec( \F_p). $$Let $r\geq 1$ be an integer. Consider the projective bundle \[f: \P(V^{(r)}) \to  X.\]
Then, the following holds.

\begin{enumerate}
    \item {Let $R$ be  a ring where $p \neq 0 \in R$ and $p \notin R^\times$.  Then, the relative scheme   $$\P(V^{(r)})\times_{\Spec(\F_p)} \Spec(R/p)  \to \Spec(R/p) $$ does not lift to a scheme flat over $\Spec(R)$.} \item{Let $R$ be a $p$-adically complete local ring with residue field $\F_p$.\\If $p \neq 0 \in R$, then  $\P(V^{(r)})$ does not lift to a scheme flat over $\Spec(R)$.}
\end{enumerate}

\end{thm}
\begin{dem}
We prove item (1). Set $Y:=\P(V^{(r)}) $. Without loss of generality,  one can assume that $R$ is of finite-type over $\Z$, hence Noetherian. One can then also assume that $\Spec(R)$ is connected. Then $p \notin p^2R$, for  otherwise $pR \subset R$ would be an idempotent ideal, hence generated by an idempotent element by Nakayama's Lemma, contradicting the connectedness of $\Spec(R)$. Replacing $R$ with $R/p^2$, one may assume $p^2=0 \in R$.\\
The obstruction to lifting $Y$ then lies in the coherent cohomology group \[ H^2(Y\times_{ \Spec( \F_p)}\Spec(R/p)  , T_{Y/ \F_p} \otimes_{\F_p} pR),\] where  $T_{Y/ \F_p}$ stands for the tangent bundle of $Y \to \Spec(\F_p)$. Using that the projection \[Y \times_{\Spec(\F_p)} \Spec(R/p) \to Y \]  is affine, this group is identified to\[H^2(Y , T_{Y/ \F_p} \otimes_{\F_p} pR).\] Since the arrow of $\F_p$-vector spaces 
$$\fonctionnoname{\F_p}{pR}{1}{p}$$
is a (split) injection by assumption,  liftability over $R$ implies liftability over $\Z/p^2$. Thenceforward,  assume $R=\Z/p^2$.\\ Assume that $Y_2$ is a lift of $Y$, flat over $\Z/p^2$. Over $Y$, there is the tautological extension of vector bundles \[\mathcal E: 0 \to \mathcal H \to  f^*(V^{(r)})  \to \mathcal O(1) \to 0.\]
The obstruction to lifting $ \mathcal H$ to a vector bundle over $Y_2$ lies in $$H^2(Y,\End( \mathcal H)) =\Ext^2_Y(\mathcal H,\mathcal H).$$ This cohomology group vanishes. To see why, use Leray's spectral sequence and the last two items of Lemma \ref{calcul3} here below.  Thus,  $ \mathcal H$  lifts to a vector bundle  $ \mathcal H_2$  over $Y_2$. Similarly, since $H^2(Y,\mathcal O_Y)=0$, the line bundle  $\mathcal O(1)$ lifts to a  line bundle $\mathcal O_2(1)$ over $Y_2$ (in fact unique up to isomorphism, since $H^1(Y,\mathcal O_Y)=0$).\vspace{0.3cm}\\
Claim: $\mathcal E$ lifts to an extension  of vector bundles over $Y_2$, \[\mathcal E_2: 0 \to \mathcal H_2 \to  V^{[r]}_2  \to \mathcal O_2(1)\to 0. \] Indeed, the obstruction to its existence lies in $$H^2(Y,\mathcal H (-1)) =\Ext^2_Y(\mathcal O(1),\mathcal H),$$ and  items (2) and (3) of Lemma \ref{calcul3} imply that this group vanishes. This proves the existence of $\mathcal E_2$. Then $V^{[r]}_2$ is a lift of $f
^*(V^{(r)})$, to a vector bundle over $Y_2$. Using Proposition \ref{exoconnection}, since $r \geq 1$, one sees that $f_*(V^{(r-1)})$ lifts to a $\W_2$-bundle over $Y$. Denoting by $W^{[r-1]}_2$ such a lift, one gets an extension of $\W_2$-modules over $Y$, \[0 \to \frob_*(f
^*(V^{r})) \to W^{[r-1]}_2  \to  f^*(V^{(r-1)}) \to 0, \]  having $\kappa= \frob$. Applying $f_*$ to this extension, using the projection formula and $R^1f_*(\mathcal O_Y)=0$, one gets the extension of $\W_2$-modules over $X$ \[0 \to \frob_*(V^{r}) \to f_*(W^{[r-1]}_2)  \to  V^{(r-1)} \to 0,  \] having $\kappa= \frob$ as well. Thus, $f_*(W^{[r-1]}_2)$ is a lift of $V^{(r-1)}$, to a $\W_2$-bundle over $X$. This contradicts Theorem \ref{nonlift}.\\
To finish the proof, it remains to prove that (1) implies (2). It suffices to prove that $Y$ is a rigid $\F_p$-variety (i.e. that it has no nontrivial deformations), for then any deformation of $Y$ over $R/p$ has to be trivial. One needs to show that $H^1(Y,T_{Y/\F_p})=0.$ The first fundamental sequence for $f:Y \to X$ (a smooth morphism of smooth $\F_p$-varieties) provides an exact sequence of vector bundles over $Y$ $$ 0 \to T_{Y/X} \to T_{Y/\F_p} \to f^*(T_{X/\F_p}) \to 0.$$ By d\'evissage, it is enough to show that $ H^1(Y,f^*(T_{X/\F_p}))$ and $H^1(Y,T_{Y/X})$  both vanish. Since $f$ is a projective bundle,   $f_*(\mathcal O_Y)=\mathcal O_X$ and  $R^1f_*(\mathcal O_Y)=0$. Using the projection formula and Leray's spectral sequence, vanishing of the first group boils down to that of  $ H^1(X,T_{X/\F_p})$. Recall that the tautological extension of vector bundles over $X$, $$ 0 \to V \to \mathcal O_X^n \to W \to 0,$$  gives rise to a natural iso $$T_{X/\F_p}=\Hom_{\mathcal O_X}(V,W).$$  The vanishing of   $ H^1(X,\Hom_{\mathcal O_X}(V,W))$ then follows from item (1) of Lemma \ref{calcul4}. To prove vanishing of  $H^1(Y,T_{Y/X})$, it suffices to prove that of $H^1(X,f^*(T_{Y/X}))$ and $R^1f^*(T_{Y/X})$. Since $f$ is the projective bundle of $V^{(r)}$, item (4) of Lemma \ref{calcul3} provides an exact sequence $$0 \to \mathcal O_X \xrightarrow{\lambda \mapsto  \lambda \Id} \End(V^{(r)}) \to f^*(T_{Y/X}) \to 0.$$ Vanishing of $H^1(X,f^*(T_{Y/X}))$  follows by d\'evissage, using that of $H^2(X,\mathcal O_X)$ (classical) and $H^1(X,\End(V^{(r)}))$ (item (5) of Lemma \ref{calcul2}). Vanishing of $R^1f^*(T_{Y/X})$ is item (5) of Lemma \ref{calcul3}, for $i=1$.
\end{dem}

\section{Some computations in coherent cohomology of flag schemes.}\label{sec computations}
Let us deal with the (more or less classical) cohomological computations used in this paper. Let $V$ be a vector bundle of rank $n\geq 2$, over a scheme $S$. \\In many  applications, $S$ is just the spectrum of a field of characteristic $p$.

\subsection{Flag schemes and classical cohomological tools.}
\begin{defi}[Flag schemes]\label{defiflag}\hfill\\
 Let \[ 1 \leq n_1 < \ldots < n_s < n \] be a strictly increasing sequence of integers. Denote by $$F(=F_{n_1,\ldots ,n_s }): \Fl(n_1,\ldots, n_s, V) \to S$$ the  scheme of flags of sub-bundles of $V$, of dimensions $n_1,\ldots,n_s$.   Denote by
\[0 \subset  \mathcal V_{n_1} \subset  \ldots \subset \mathcal V_{n_s} \subset \mathcal V_n= F^*(V) \]  the tautological flag  over $\Fl(n_1,\ldots, n_s, V)$.

Denote $\Fl(1,2,\ldots,n-1,V)$ simply by $\Fl(V)$; it is the scheme of complete flags of the vector bundle $V$. For $1 \leq i \leq j \leq n$, denote by \[\mathcal V_{j/i}:= \mathcal V_j / \mathcal V_i \]and  \[\mathcal L_{i}:= \mathcal V_{i} / \mathcal V_{i-1} \] the natural quotients. \\
 For an arbitrary sequence of relative integers $a_1, \ldots, a_n$,  put \[\O(a_1, \ldots, a_n):= \mathcal L_{1} ^{\otimes a_1}  \otimes \ldots \otimes \mathcal L_n ^{\otimes a_n}.  \]

\end{defi}

\begin{prop}\label{CohFlag}
  Denote by $$f: \P(V) \to S$$ the projective bundle of $V$,  by $\O(1)$ its twisting sheaf, and $$F:\Fl(V) \to S$$ its complete flag scheme. (Note that $\O(1)=\mathcal L_n$, with the notation above.) \\ Let $m\geq 0$  and   $a=(a_1, \ldots, a_n) \in \Z^n$ be   integers. The following holds. 
\begin{enumerate}
    \item {$f_*( \O(m))=\Sym^m(V)$, and $R^if_*( \O(m))= 0$, for $i \geq 1$.} \item{$R^if_*( \O(-m))= 0$ for $i < n-1$.} \item{$R^{n-1}f_*( \O(-m))= 0$, for $0 \leq m \leq n-1$.} \item{ $R^{n-1}f_*( \O(-m))= \Gamma^{m-n}(V^\vee)  \otimes \Det(V)^\vee$, for $m \geq n$.}  \item{Computations similar to (1)-(4),  but over $\Fl(V)$ in place of $\P(V)$. Just replace $\O(*)$ by $\O(0,\ldots,0,*)$, or dually $\O(-*,0,\ldots,0)$. When $n=2$,  one gets  a natural complete flag on the vector bundle $$F^*(R^1F_*( \O(m,0))).$$ Its  graded pieces are the line bundles $\O(n_1,n_2)$,  $n_1,n_2 \geq 1$,  $n_1+n_2=m$.} \item{If $a$ is  not an increasing sequence, then \[ F_*(\O(a_1, \ldots, a_n))=0. \]}  \item{If $a$ is an increasing sequence, then, for all $i \geq 1$, \[ R^i(F_*)(\O(a_1, \ldots, a_n))=0. \]}
\end{enumerate}

\end{prop}
\begin{dem}
For items (1) to (4), see \cite[Tag 30.8]{SP}. Item (5) follows, using Leray's spectral sequence and the projection formula, w.r.t. the factorisation $$\Fl(V) \to \P(V) \to S.$$
For items (6) and (7), see  \cite[Proposition 1.4.5]{B}.
\end{dem}
\subsection{Some computations.}
\begin{lem}\label{calcul3}
 Consider the projective bundle \[f:\P(V) \to  S\,\] with  tautological extension \[\mathcal E: 0 \to \mathcal V_{n-1} \to  \mathcal V_n(=f^*(V))  \to \mathcal L_n(=\mathcal O(1)) \to 0.\]  The following equalities hold. \begin{enumerate}
    \item {$R^if_*(\mathcal V_{n-1} )=0$, for all $i \geq 0$.}\\
    \item{$R^if_*(\mathcal V_{n-1} (-1))=0$, for all $i\neq 1$.}\\
    \item {$R^1f_*(\mathcal V_{n-1}(-1))=\mathcal O_S,$ with generator given by $\mathcal E$.}\\ \item {$f_*(\mathcal V_{n-1}^\vee(1))=\End(V)/\mathcal O_S \Id$.}\\ \item {$R^if_*(\mathcal V_{n-1}^\vee(1))=0,$ for all $i\geq 1$.}\\
    \item {$R^if_*(\End(\mathcal V_{n-1}))=0$, for $i \geq 1$.}\\
    \item { $f_*(\End(\mathcal V_{n-1}))=\mathcal O_S,$ with generator given by $\Id_{\mathcal V_{n-1}}$.} 
\end{enumerate}


    
    
   
\end{lem}
\begin{dem}
Applying $f_*$  to the surjection of vector bundles over ${\P(V)}$ $$ f^*(V)  \to \mathcal O(1), $$ one gets $$\Id: V \to V.$$
Item (1) follows, applying $f_*$ to $\mathcal E$, together with the projection formula, and the well-known equalities  $$f_*( \mathcal O_{\P(V)})=\mathcal O_S,$$  and  $$R^i f_*( \mathcal O_{\P(V)})=0,$$ for all $i \geq 1$. \\
To get items (2) to (5), consider the extensions \[\mathcal E(-1): 0 \to \mathcal V_{n-1}(-1) \to  f^*(V)(-1)  \to \mathcal O_{\P(V)} \to 0\] and \[\mathcal E^\vee(1): 0 \to  \mathcal O_Y \to  f^*(V^\vee)(1)  \to  \mathcal V_{n-1}^\vee(1) \to 0.\]  Apply $f_*(.)$ to these, using the projection formula, and the vanishing of  $$f_*( \mathcal O(-1)),R^i f_*( \mathcal O(1)),R^i f_*( \mathcal O(-1)),$$which holds  for all $i \geq 1$.
Introduce the extension  \[\mathcal E^\vee \otimes \mathcal V_{n-1}: 0 \to \mathcal V_{n-1}(-1) \to  f^*(V^\vee) \otimes \mathcal V_{n-1} \to \End(\mathcal V_{n-1}) \to 0.\]  Applying the same technique as above,  one gets items (6) and (7).\\

\end{dem}

One can perform similar computations on Grassmannians and other flag schemes. \\ If $S$ is an $\F_p$-scheme, some items above remain valid for Frobenius twists- depending  on the cohomological degree and  on the dimension of the tautological bundles.\\ Positive examples are treated in a  systematic manner in  the Lemmas below.

\begin{lem}\label{calcul4}
Assume that $S$ is an $\F_p$-scheme.\\
 For $1 \leq m  \leq  n-2$, consider  the Grassmannian $$f: \Gr(m,V) \to S,$$ over which there is the tautological extension of vector bundles \[\mathcal E: 0 \to \mathcal V_m \to f^*(V) \to \mathcal V_{n/m} \to 0. \]

    The following is true. \begin{enumerate}
        \item{For all $i\geq 1$, $R^i f_*(\mathcal V_m^\vee \otimes \mathcal V_{n/m})=0$.}  \item{For all $r\geq 0$, $\Ext^1_{\mathcal O_{\Gr(m,V)}}(\mathcal V_m^{(r)},\Sym^{p^r}(\mathcal V_m))=0.$}  \item{ For all  $i \geq 0$, $f_*(\End(  \Gamma^i(\mathcal V_m) ))=\O_S \Id.$}  \item{ For all $r\geq 1$, $\Hom_{\mathcal O_{\Gr(m,V)}}(\mathcal V_m^{(r)},  \Gamma^{p^r}(\mathcal V_m) )=0.$}
    \end{enumerate}

Moreover, items (2), (3) and (4) hold, replacing $\mathcal V_m$ by $\mathcal V_m^{(s)}$, for all $s \geq 1$.
\end{lem}
\begin{dem}
Let us prove item (1).
Consider the factorisation  $$ \Fl(V) \xrightarrow{g} \Gr(m,V) \to S,$$  where $g$ is the composite of the complete flag schemes of $\mathcal V_m$ and $\mathcal V_{n/m}$. One has $g_*(\mathcal O_{\Fl(V)})=\mathcal O_{\Gr(m,V)}$, and   $R^i g_*(\mathcal O_{\Fl(V)})=0$ for all $i\geq 1$. This is easily derived from the similar well-known formulas for projective bundles. \\Hence, using the projection formula, it suffices to prove (1) after base-change to $\Fl(V)$ (that is to say, with $F$ in place of $f$). There, the  vector bundle  $\mathcal V_m^\vee \otimes \mathcal V_{n/m}$ acquires a complete filtration, whose graded pieces are line bundles of the shape $\O(0,\ldots,0,-1, 0,\ldots,0,1,0,\ldots,0)$. These have $R^i F_*(.)=0$, for all $i \geq 1$. Indeed, for the line bundle $\O(-1, 0,\ldots,0,1)$, this is item (7) of Lemma \ref{CohFlag}. In all other cases, the sequence of integers inside  $\O(.)$ contains two consecutive terms, which are either $(0,-1)$ or $(1,0)$. Say their indices are $j,j+1$, so that the corresponding line bundles on $\Fl(V)$ are $\mathcal L_j$ and  $\mathcal L_{j+1}$ Consider the factorisation $$ \Fl(V) \xrightarrow{h} \Fl(1,\ldots, j-1, j+1,j+2, \ldots, n-1,V) \to S.$$  Then  $h$ is the $\P^1$-bundle $\P(\mathcal V_{j+1/j-1})$, with twisting sheaf $\mathcal O(1):=\mathcal L_{j+1}$ . Recall that, for all $i \geq 0$, one has $$R^i h_*(\mathcal L_{j+1}^{-1})=R^i h_*(\mathcal O(-1))=0,$$ and similarly $$R^i h_*(\mathcal L_j)=R^i h_*(\Det(\mathcal V_{j+1/j-1})(-1))=0,$$ by the projection formula. Using Leray's spectral sequence for the factorisation above, one indeed concludes that $R^i F_*(.)=0$, for all graded pieces and for all $i \geq 1$.  The claim follows by d\'evissage on the filtration. Let us prove items (2), (3) and (4) as stated (i.e. when $s=0$). The proof for $s \geq 1$ is the same. For (2), arguing as above, one gets that $$\Ext^1_{\mathcal O_{\Gr(m,V)}}(\mathcal V_m^{(r)},\Sym^{p^r}(\mathcal V_m))=\Ext^1_{\O_{\Fl(V)}}(\mathcal V_m^{(r)},\Sym^{p^r}(\mathcal V_m)).$$ Using the projection formula and item (1) of Proposition \ref{CohFlag}, one then sees that $$\Ext^1_{\O_{\Fl(V)}}(\mathcal V_m^{(r)},\Sym^{p^r}(\mathcal V_m))= \Ext^1_{\O_{\Fl(V)}}(\mathcal V_m^{(r)},\mathcal L_m^{p^r}).$$ By d\'evissage along the natural filtration of $\mathcal V_m$, it  remains to prove vanishing of $$\Ext^1_{\O_{\Fl(V)}}(\mathcal L_i^{p^r},\mathcal L_m^{p^r}))=H^1(\Fl(V),\mathcal L_i^{-p^r} \otimes \mathcal L_m^{p^r}),$$ for $i=1,\ldots, m$. This is clear if $i=m$. If $i < m$, consider the factorisation $$ \Fl(V) \to \Fl(1,\ldots,m-1, m,V)  \xrightarrow{h} \Fl(1,\ldots, m-1,V)\to S,$$ where $h$ is the projective bundle of the vector bundle dual to $\mathcal V_n /\mathcal V_{m-1}$. Its  twisting sheaf is $\O(1):=\mathcal L_m^{-1}$. Since $m \leq n-2$,  one has $$R^jh_*(  \mathcal L_m^{p^r})=R^jh_*( \O(-p^r))= 0,$$ for $j=0,1$, by item (2) of Proposition \ref{CohFlag}. Once more, we conclude using Leray spectral sequence and projection formula. By duality (between $\Sym$ and $\Gamma$), the equality in item (3) is equivalent to $$f_*(\End(  \Sym^i(\mathcal V_m) ))=\O_S \Id,$$ which we now prove. As above, this first reduces to $$f_*(\End_{\Fl(V)}(  \Sym^i(\mathcal V_m) ))=\O_S \Id,$$  then to $$f_*( \Sym^i(\mathcal V_m)^\vee \otimes \mathcal L_m^i )=\O_S \pi_m ^i,$$ where $\pi_m: \mathcal V_m \to \mathcal L_m$ is the natural surjection. The vector bundle $\Sym^i(\mathcal V_m)$ has a natural complete filtration, with graded pieces the line bundles \[\O(a_1, \ldots, a_m,0,\ldots,0),\] one for each partition $a_1+ \ldots+a_{m-1}+ a_m=i$. By d\'evissage, it suffices to prove $$f_*( \O_{\Fl(V)})=\O_S,$$ for the partition $0+ \ldots+0+i=i$, and $$f_*( \O(-a_1, \ldots,\ldots,-a_{m-1},i-a_m,0 \ldots,0))=0,$$ for all other partitions. The former  is clear. For the latter, observe that $i-a_m <0$, and use item (6) of Proposition \ref{CohFlag}. \\To prove item (4), pick $$\phi \in \Hom_{\mathcal O_{\Gr(m,V)}}(\mathcal V_m^{(r)},  \Gamma^{p^r}(\mathcal V_m) ).$$ Consider the composite $$\psi: \Gamma^{p^r}(\mathcal V_m)  \xrightarrow{\frob_{\mathcal V_m}^r}\mathcal V_m^{(r)}   \xrightarrow{\phi} \Gamma^{p^r}(\mathcal V_m).$$ By item (3), one gets $\psi= \lambda \Id$, where  $$\lambda \in H^0(S,\O_S) = H^0(\Gr(m,V),\O_{\Gr(m,V)}).$$   It suffices to prove $\lambda=0$, for then the surjectivity of $\frob_{\mathcal V_m}^r$ implies  $\phi=0$. To finish, observe  that  vanishing of $\lambda$ can be checked over each affine open $U \subset \Gr(m,V)$, where the vector bundle $\mathcal V_m$ is trivial. Over such a $U=\Spec(R)$, choosing an $R$-basis of $\mathcal V_m$ identifies the map of $R$-modules  $$\frob_{\mathcal V_m}^r: \Gamma^{p^r}(\mathcal V_m)  \to \mathcal V_m^{(r)}$$ to a projection $$\pr: R^M \xrightarrow{(x_1,\ldots,x_M) \to (x_1,\ldots,x_m)} R^m.$$  It then becomes obvious that $\lambda \Id_{R^M}$  factors through $\pr$,  if and only if $\lambda=0$.
\end{dem}

\begin{lem}\label{calcul2}
Assume that $S$ is an $\F_p$-scheme. Let $1 \leq m \leq n-2$ be an integer.\\
Denote by $$F: \Fl(V) \to S$$ the complete flag scheme of $V$, and by \[ \pi_m: \mathcal V_m \to \mathcal L_m\] the natural surjection of vector bundles over $\Fl(V)$, with kernel $\mathcal V_{m-1}$. \\Denote by \[\mathcal E_{m-1,m} : 0 \to \mathcal V_{m-1} \to \mathcal V_m \xrightarrow{\pi_m} \mathcal L_m \to 0\] the tautological extension of vector bundles over  $\Fl(V)$.\vspace{0.3cm}\\
Let $r\geq 0$ be an integer. The following equalities hold.
\begin{enumerate}
    \item $F_*(\mathcal V_m^{(r)\vee} \otimes\mathcal L_m^{(r)})=\O_S \pi_m^{(r)}.$
    \item $R^1F_*(\mathcal V_m^{(r)\vee }\otimes \mathcal L_m^{(r)})=0.$  
    \item $R^1F_*(\mathcal L_m^{(r)\vee }\otimes \mathcal V_{m-1}^{(r)})=\O_S \mathcal E_{m-1,m}^{(r)} .$
    \item $F_*(\mathcal L_m^{(r)\vee}\otimes \mathcal V_{m-1}^{(r)})=0 .$
    \item $R^1F_*(\End(\mathcal V_m^{(r)}))=0.$
\end{enumerate}
\end{lem}

\begin{dem} Proofs of items (1) and (2), are the same as that of items (2) and (3) of Lemma \ref{calcul4}. 
 For $m=2$, using Leray's spectral sequence w.r.t. \[ \Fl(V) \stackrel  f \to \Fl(3,4,\ldots, n-1, V) \to S,\] one deduces (3) and (4)  from Lemma \ref{LemExt} here below (applied to  the vector bundle $\mathcal V_3$, over $\Fl(3,4,\ldots, n-1, V)$). For $m = 3$, consider  (the long exact sequence in cohomology obtained by applying $f_*(.)$ to) the natural extension $$ 0 \to \mathcal L_3^{(r)\vee }\otimes \mathcal V_2^{(r)} \to \mathcal L_3^{(r)\vee }\otimes \mathcal V_3^{(r)} \to  \mathcal O_{\Fl(V)} \to 0. $$ By item (2)  of Proposition \ref{CohFlag}, and the projection formula,  one gets,  for $i=0,1$,  $$R^i f_*(\mathcal L_3^{(r)\vee }\otimes \mathcal V_3^{(r)} )= R^i f_*(\mathcal L_3^{-p^r })\otimes \mathcal V_3^{(r)}=0 .$$ Conclude using Leray's spectral sequence. The proof for $m \geq 4$ is similar.
 
 Let us prove item (5). The case $m=1$ boils down to vanishing of $R^1F_*(\mathcal O_{\Fl(V)})$.  The case $m \geq 2$ follows by induction, like this. Consider the exact sequence $$ 0 \to  \mathcal V_m^{(r)\vee }\otimes \mathcal V_{m-1}^{(r)} \to \End(\mathcal V_m^{(r)}) \to    \mathcal V_m^{(r)\vee }\otimes \mathcal L_m^{(r)} \to 0. $$ Using item (2), it remains to prove vanishing of $R^1F_*(\mathcal V_m^{(r)\vee }\otimes \mathcal V_{m-1}^{(r)})$.\\ To do so, consider the natural extension  $$ 0 \to  \mathcal L_m^{-p^r } \otimes \mathcal V_{m-1}^{(r)} \to \mathcal V_m^{(r)\vee }\otimes \mathcal V_{m-1}^{(r)}\to  \End(\mathcal V_{m-1}^{(r)}) \to 0. $$ Using the induction hypothesis, it remains to prove vanishing of the induced arrow  $$ h: R^1F_*(\mathcal L_m^{-p^r } \otimes \mathcal V_{m-1}^{(r)}) \to  R^1F_*(\mathcal V_m^{(r)\vee }\otimes \mathcal V_{m-1}^{(r)}).$$ Using item (3), the left side is a free $\mathcal O_S$-module of rank one, whose generator $ \mathcal E_{m-1,m}^{(r)} $ is tautologically killed by $h$. The proof is over.
\end{dem}

The next Lemma is a rigorous formulation of the 
\begin{mot}
   A ``sufficiently general'' extension of line bundles has no non-trivial automorphism, and  its  class generates the cohomology group in which it dwells.
 
\end{mot}
 \begin{lem}\label{LemExt}
Assume that $n=\dim(V)=3$, and that $S$ is an $\F_p$-scheme.\\ Consider  the composite \[f:\Fl(V)= \Fl(1,2,V) \stackrel {f_1} \to  \Fl(2,V)=\P(V) \stackrel {f_2} \to S.\] Then $$ f_*(\mathcal O(p^r,-p^r,0))=0,$$ and $$ R^1f_*(\mathcal O(p^r,-p^r,0))=\mathcal O_S,$$ with generator given by the class of the tautological extension  \[\mathcal E_{1,2}^{(r)}: 0 \to \mathcal L_1^{(r)} \to \mathcal V_2^{(r)} \to  \mathcal L_2^{(r)}  \to 0. \]
 \end{lem}
\begin{dem}
The usual computation of cohomology of $\P^1$-bundles, applied to $f_1$, gives $$ (f_1)_*(\mathcal O(p^r,-p^r,0))=0,$$ and $$ R^1(f_1)_*(\mathcal O(p^r,-p^r,0))=\Gamma^{2p^r-2}(\mathcal V_2) \otimes \Det (\mathcal V_2)^{1-p^r}.$$ Indeed, the arrow $f_1$ is a $\P^1$-bundle, the projective bundle of $\mathcal V_2$, and  $$\mathcal O(p^r,-p^r,0)=\mathcal O(-2p^r)\otimes \Det(\mathcal V_2)^{\otimes p^r},$$ where $\mathcal O(1)$ denotes the usual twisting sheaf of a  $\P^1$-bundle. Over $\Fl(V)$, the vector bundle $\Gamma^{2p^r-2}(\mathcal V_2) \otimes \Det (\mathcal V_2)^{1-p^r}$ has a natural filtration, with graded pieces the line bundles $\mathcal O(-a,a,0)$, where $1- p^r \leq a \leq  p^r -1$. These have $f_*(.)=0$, except  when $a=0$, for then $f_*(\mathcal O_{\Fl(V)})=\mathcal O_S$. Using (a tiny portion of) Leray's spectral sequence for $f=f_2 \circ f_1$,  we get an injective arrow of coherent $\mathcal O_S$-modules \[ R^1f_*(\mathcal O(p^r,-p^r,0)) \stackrel \sim \to (f_2)_*(\Gamma^{2p^r-2}(\mathcal V_2) \otimes \Det (\mathcal V_2)^{1-p^r}) \hookrightarrow \mathcal O_S.\] It admits the natural arrow 
$$\fonctionnoname{\mathcal O_S}{ R^1f_*(\mathcal O(p^r,-p^r,0))}{1}{<\mathcal E_{1,2}^{(r)}>},$$
as a splitting, which is hence an isomorphism.
\end{dem}

\begin{rem}
Assumptions and notation being those of Lemma \ref{LemExt}, one can show that $ R^1f_*(\mathcal O(a,-a,0))=0,$ for every integer $a\geq 0,$ which is not a $p$-th power.
\end{rem}


\begin{thebibliography}{ABCD}

\bibitem[BGA]{BGA} A. Bertapelle, C. Gonz\'alez-Aviles, \textit{The Greenberg functor revisited}, Eur. J. Math. 4 (2018), 1340--1389.

\bibitem[Bo]{BO} J. Borger, \textit{The basic geometry of Witt vectors, I: The affine case}, Algebra Number Theory 5 (2011), 231--285.

\bibitem[B]{B} M. Brion, \textit{Lectures on the Geometry of Flag Varieties}, in Trends in Math.: Topics in Cohomological Studies of Algebraic Varieties,  Birkh\"auser, 33--85.

\bibitem[DK]{DK} C. Davis, K. Kedlaya, \textit{On the Witt vector Frobenius}, Proc. Amer. Math. Soc. 142(7) (2014), 2211--2226.

\bibitem[DI]{DI} P. Deligne, L. Illusie, \textit{Rel\`evements modulo $p^2$ et d\'ecomposition du complexe de de Rham}, Inv. Math. 89 (1987), 247--270.

\bibitem[Fe]{Fe} D. Ferrand, \textit{Un foncteur norme}, Bull. Soc. Math. France 126 (1998), 1--49.

\bibitem[Ha]{Ha} M. Hazewinkel, \textit{Witt vectors}, in Handbook of algebra, Vol. 6, 2009, 319--472.

\bibitem[Il]{Il} L. Illusie, \textit{Complexe Cotangent et D\'eformations I}, Lecture Notes in Math. 239, 1971.

\bibitem[Le]{Le} H. Lenstra, \textit{Construction of the ring of Witt vectors}. 2022. Available at:\\
\url{http://www.math.leidenuniv.nl/~hwl/}

\bibitem[MR]{MR} V. B. Mehta, A. Ramanathan, \textit{Frobenius splitting and cohomology vanishing of Schubert varieties}, Annals of Math. 122 (1985), 27--40.
 
\bibitem[MS]{MS} V. B. Mehta, V. Srinivas, \textit{Varieties in positive characteristic with trivial tangent bundle}, Comp. Math. 64 (1984), 191--212.

\bibitem[Se1]{Se1} J.-P. Serre, \textit{Corps locaux}, 2nd edition, Hermann, Paris, 1968.

\bibitem[Se2]{Se2} J.-P. Serre, \textit{Sur la topologie des vari\'et\'es alg\'ebriques en caract\'eristique $p$}, Symposium de topologie alg\'ebrique, Mexico, 1956, 24--53.

\bibitem[Sch]{Sch} S. Schr\"oer, \textit{The Deligne-Illusie theorem and exceptional Enriques surfaces}, Eur. J. Math. 7 (2021), 489--525.

\bibitem[SP]{SP} \textit{The Stacks Project}, https://stacks.math.columbia.edu/

\bibitem[TO]{TO} T. Sekiguchi, F. Oort, \textit{On the deformation of Witt groups to tori}, Algebraic and topological theories, Kinokuniya Co., Tokyo, Japan, 1986, 283--298.

\bibitem[Yo]{Yo} F. Yobuko, \textit{Quasi-Frobenius splitting and Lifting of Calabi-Yau varieties in characteristic $p$}, Math. Z. 292 (2019), 307--316.

\bibitem[Z]{Z} M. Zdanowicz, \textit{Arithmetically rigid schemes via deformation theory of equivariant vector bundles},  Math. Z. 297 (2021), 361--387.
\end{thebibliography}
\end{document}